\documentclass{article}
\usepackage[a4paper]{geometry}
\usepackage[T1]{fontenc}
\usepackage{lmodern}

\usepackage{graphicx, caption, subcaption, float} % Required for inserting images
\usepackage{amsmath, amssymb, amsfonts}
\usepackage{hyperref, booktabs}
\usepackage{empheq, mdframed}
\newmdtheoremenv{theo}{Theorem}
\usepackage{physics}
\graphicspath{ {./images_png/} }

\numberwithin{equation}{section}

\usepackage[
backend=biber,
style=phys,
sorting=none,
minnames=2,
maxnames=2,
biblabel=brackets, % APS
pageranges=false,
chaptertitle=false
]{biblatex}
\DeclareFieldFormat{url}{\ifentrytype{online}{URL: \href{#1}{#1}}{}}

\usepackage{tikz}
\usepackage{dirtytalk}
\usepackage{appendix}

\usepackage{authblk}

\title{Deducing the Talbot Effect from Electrodynamics}
\author[]{Gabriel Ybarra Marcaida}
\affil[]{Email: \texttt{gabriybarra@gmail.com}}
\date{}

\begin{document}

\maketitle

\section*{Abstract}

We propose a model based on Maxwell's equations to describe the Talbot effect. After solving the model analytically, we prove that its solution is equivalent to the one amply found in the literature in the asymptotic limit. By considering the paraxial limit, we showcase how the rational Talbot effect arises, and derive an ideal paraxial distribution that would give such effects. We give a proof of $L^2$-convergence to these ideal distributions for any such grating. The universal nature of the Talbot effect makes our findings pertinent for a variety of physical systems and technological applications.

\newpage

\section{Introduction}

Lord Talbot discovered in 1836 that in the near field of a diffraction grating a (Talbot) plane depicting a clear replica of the diffraction grating appeared without the need of using lenses, a phenomena that occurred at periodic distances away from the grating \cite{Talbot1836}. Half a century later, Lord Rayleigh showed in \cite{Rayleigh1881} that this periodic distance, which he baptised as the \emph{Talbot distance}, could be computed through the relation $z_T = 2d^2/\lambda$, where $\lambda$ is the wavelength of the source of light and $d$ is the distance between gratings.\smallskip

It was later found that self-images with a smaller period --or intermediary Talbot patterns-- appear at every rational multiple of the Talbot distance. This is now called the \emph{fractional Talbot effect} \cite{cowley_fourier_1957, Hiedemann1959, Winthrop1965}. At any distance $z= \frac{p}{q} \, \frac{z_T}{2}$, where both $p$ and $q$ are coprime integers, $q$ well-defined subimages will be formed. $q$ is called the \emph{order} of the Talbot subimage.\smallskip

Further progress was made in 1996 in a paper by Berry and Klein \cite{Berry1996}. First, they showed that the phases of the Talbot images are related to the quadratic Gauss sums, a number-theoretical object. Secondly, and maybe more interestingly, they proved that the gratings with sharp-edge slits produce fractal intensity graphs at irrational multiples of the Talbot distance.\smallskip

Moreover, an effect analogous to Talbot's can be found in other media like Bose-Einstein condensates \cite{Deng1999, Ryu2006}, spin-waves \cite{Golebiewski2020}, water \cite{Bakman2018, Sungar2018, Rozenman2022} or acoustic waves \cite{Saiga85}, among others. In particular, a \emph{temporal} Talbot effect arises in quantum mechanics, where it is known as the \emph{quantum revival of the wave function}. The rational Talbot effect can be used to construct a universal set of quantum gates and test the $D$-dimensional Bell inequalities, as was shown in \cite{Jimenez2015, Barros2017}. Another use of the fractional Talbot effect is that it allows for the factorisation of integers, as explained in \cite{Clauser1996}. More recently, it was shown in \cite{Chamizo2024} that the Dirac equation exhibits exact quantum revivals, which might find application in solid state physics.\smallskip

Similar effects can also be found in systems described by the non-linear Schrödinger equation. This equation connects the Talbot effect to other phenomena like rogue waves, which are spontaneous and huge surface waves \cite{Nikolic2019}, or the evolution of vortex filaments \cite{DeLaHoz2014}. Simulations of the evolution of the vortex filaments, where we can see the self-imaging phenomenon, are available at \cite{KumarVideos}. All of this shows that our findings are not only relevant for diffraction optics, but for a wide variety of research areas and technological applications.\smallskip

The theoretical approach to the Talbot effect is rather similar in all the papers cited above. Indeed, a scalar theory is used, whose solution is given by solving the Helmholtz equation. This solution is (see \cite{Mejias1990, Berry1996})
\begin{equation}\label{eq:stationary_solution_cited}
    U (x,z) = \sum_n \hat{g}_n \cos{k_n x}\times \begin{cases}
        e^{- iz\sqrt{\omega^2 - k_n^2} }, &  0< k_n \leq\omega \\ 
        e^{-z\sqrt{k_n^2 - \omega^2}}, &  \omega< k_n <\infty,
    \end{cases}
\end{equation}
where $k_n = \dfrac{2\pi}{d}n$ and $\hat{g}_n$ is a Fourier coefficient that depends on the grating. This agrees with the solution obtained by instead considering the Fresnel integral, as was proven in \cite{Matsutani2003}. The problem with this approach is that it relies on solving the Helmholtz equation with Kirchhoff conditions, which means that we impose the value of the field and its derivative on the screen. But Poincaré showed that imposing these two conditions simultaneously leads to a wave that vanishes everywhere \cite{Buchwald2006, Poincare1887}! Alternatively, \emph{Sommerfeld radiation condition} must be used, which states that incoming waves are unphysical and must be discarded \cite{sommerfeld_greensche_1912}. Although being physically sensible, this approach lacks mathematical justification.\smallskip

Our contributions to the topic are three. First, we propose a time-dependent model based on classical electrodynamics, which we solve exactly and which allows us to identify the field that is described by \eqref{eq:stationary_solution_cited} as the electromagnetic four-potential $A^y$. Critically, this model does not need us to impose Kirchhoff's boundary conditions nor Sommerfeld's radiation condition. Next, we prove that our solution agrees with the well-known \eqref{eq:stationary_solution_cited} in the limit $t \to\infty$. And, finally, we prove the $L^2$ convergence of the solution to the ideal distribution in the paraxial limit.

The paper is organised as follows. In Section \ref{sec:our_model} we give a well-posed PDE problem which describes the Talbot effect, which we justify by considering classical electrodynamics an Talbot's experimental setting. Next, we obtain in Section \ref{section:transient_solution} an analytical solution to this model that establishes the formation and evolution of the Talbot effect. In Section \ref{section:stationary_solution} we show that we can recover \eqref{eq:stationary_solution_cited} by considering the asymptotic temporal behaviour of our solution, thus giving a stronger physical and mathematical footing to Talbot's effect. In Section \ref{section:graphing_carpets} we use a numerical approach to study the solutions to the wave equation for several parameters, from which we extract what the Talbot effect looks like in the zero wavelength --or paraxial-- limit. We propose in Section \ref{section:paraxial} a distribution that shows an exact rational Talbot effect and prove that the solution to the wave equation converges in the $L^2$ sense to it in the paraxial limit for a wide class of gratings. In Section \ref{section:ideal}, we study the thin grating case, which relates to the fractal Talbot effect and derive an ideal distribution that describes it. Finally, we show that this distribution displays \emph{dark paths}, that is, paths along which there is total darkness.\smallskip

\section{Our model}\label{sec:our_model}

First of all, we propose a differential equation and a set of boundary conditions that replicate Lord Talbot's experiment's features. This should define a function $u(t,x,y,z)$ that describes the Talbot effect in the region $z \geq 0$. This effect appears when a plane wavefront of light is diffracted by an infinite set of equally spaced and identical slits. As we are dealing with electromagnetic waves, we would like our solution to be compliant with Maxwell's theory. In order to do this, we recall that we may define the electromagnetic four-potential $A^\mu = (\phi, \mathbf{A})$ from which we can recover the electric and magnetic fields through the relations
\begin{equation*}
    \mathbf{E} = -\nabla\phi - \partial_t \mathbf{A}, \quad \mathbf{B} = \nabla\cross\mathbf{A},
\end{equation*}
where we are using natural units ($c=1$). Maxwell equations for the four-potential are
\begin{equation}\label{eq:four-potential_maxwell}
    \partial_\nu (\partial^\nu A^\mu - \partial^\mu A^\nu) = J^\mu,
\end{equation}
where $J^\mu=(\rho, \mathbf{J})$ is the electromagnetic four-current. We may work in the Lorenz gauge, $\partial_\mu A^\mu = 0$, for which \eqref{eq:four-potential_maxwell} reduces to
\begin{equation*}
    \partial_\nu \partial^\nu A^\mu = J^\mu.
\end{equation*}
Now, in our case there are no electromagnetic sources, as we expect those to be somewhere in the region $z < 0$, and thus are in the vacuum. Furthermore, we consider that the grating modifies the propagation of light without introducing any source. We ignore any complications related to the physical nature of the grating, notably the fact that a real slit acts as a charged sheet, has a finite thickness and is not be perfectly conducting. A more rigorous treatment of these electromagnetic effects can be found in \cite{Noponen1993}.\smallskip

We then have $J=0$ and thus $A^\mu$ satisfies a wave equation
\begin{equation}\label{eq:four-potential_maxwell_lorenz}
    \partial_\nu \partial^\nu A^\mu \equiv \partial_{tt} A^\mu - \nabla^2 A^\mu  = 0.
\end{equation}\smallskip

We have thus found that the differential equation that $u$ satisfies must be the wave equation, which is
\begin{equation}\label{wave}
    \partial_{tt} u = \partial_{xx} u + \partial_{yy} u + \partial_{zz} u.
\end{equation}
This can be simplified by considering that we are dealing with the ideal case of infinitely long gratings, which makes $u(t,x,y,z) = u(t,x,y',z)$, for all $y,y'\in\mathbb{R}$. So the $y$-coordinate can be disposed of: $u(t,x,y,z)\equiv u(t,x,z)$. We will also assume, without loss of generality, that the diffraction grating is set along the $x$ axis, at $z=0$. These assumptions allow us to identify 
\begin{equation*}
    A^\mu(t,x,z) = \left(0, 0, u(t,x,z), 0 \right).
\end{equation*}
This way we find that $\mathbf{E} = -\partial_t A^y \mathbf{\hat{y}}$ and $\mathbf{B} = -\partial_z A^y \mathbf{\hat{x}} + \partial_x A^y \mathbf{\hat{z}}$, so that the propagation of light is in the $XZ$-plane, as expected. Furthermore, we see that $\partial_\mu A^\mu = \partial_y A^y = 0$ as required by the gauge fixing condition.\smallskip

Now, as there are infinitely many slits, the solution must be periodic, with a period $d$, equal to the distance between two slits, $u(t,x,z)=u(t,x+d,z)$. So we can restrict ourselves to studying the interval $0\leq x\leq d$. Moreover, the wave must be equal at the left and the right side of every slit, so we must find an even function on $x$. This immediately gives us the boundary condition
\begin{equation}\label{eq:condicion_borde}
    \left.\partial_x u\right\vert_{x=0} = \left.\partial_x u\right\vert_{x=d/2} = 0.
\end{equation}
Knowing this, we can restrict ourselves to the region $0 < x < d/2$, $z\geq 0$, $t\geq 0$, as the complete solution may be recovered from the one in this region using the symmetries of the problem.\smallskip

Now, we also know that the source is positioned at the grating, so the value of the function $u$ at the grating at every instant of time should be known to us. We will consider the general case where the profile of light at the grating will be described by a function $f(t,x)$, which must be periodic in $x$, with a period of $d$. Nevertheless, in order to simplify the equations, we will consider that the source of light is turned on at the instant $t=0$. We can thus write
\begin{equation}\label{eq:condicion_rendija}
    u(t,x,y,z=0) = f(t,x) := g(x)\, h(t)\, \theta(t),
\end{equation}
where $g(x)$ is periodic with period $d$, $\omega$ is the frequency of the light and $\theta(t)$ is Heaviside's theta function. In the monochromatic case we take $h(t) := \sin \omega t$.\smallskip

The other boundary condition for the $z$-direction is
\begin{equation}\label{eq:condicion_fondo}
    \lim_{z\to\infty} \abs{u(t,x,z)} = L \neq \infty ,
\end{equation}\smallskip
which just means that $u$ is bounded at infinity. This will ensure the unicity of the solution, as we will see later.\smallskip

\begin{figure}[H]
    \centering
    \begin{tikzpicture}
        \draw (0, 4) -- (-5,4) (-2,4) node[anchor=south, align=center] {At $x=d/2$:  $\left.\partial_{x} u \right\vert_{x=d/2}=0$};
        \draw (-5,4) --  (-5,0) (-5.2,2) node[anchor=south, align=center, rotate=90] {At the grating:\\ $u(t,x,0)=f(t,x)$};
        \draw (-5,0) -- (0, 0) (-2,0) node[anchor=north, align=center] {At $x=0$: $\left.\partial_{x} u \right\vert_{x=0}=0$};
        \draw (-2,2) node[align=center]{Inside the domain:\\ $\partial_{tt} u = \laplacian{u}$};
        \draw[dashed] (0,0) -- (1,0) -- (1,4) -- (0,4) (1.2,2) node[anchor=west, align=center] {At infinity:\\$u(x,z\rightarrow\infty,t)$\\is bounded};

        %\node[anchor=east] at (-5,0) {$x=0$};
        %\node[anchor=east] at (-5,4) {$x=d/2$};
        
        % Lines
        \draw [->] (-6.5,0.5) -- (-6.5,3.5);
        \node[anchor=south, rotate=90] at (-6.5,2) {$x$-direction};

        \draw [->] (-4,-1) -- (0,-1);
        \node at (-2,-1.3) {$z$-direction};
    \end{tikzpicture}
    \caption{Sketch of the differential equation and boundary conditions obeyed by the system.}\label{problem_diagram}\vspace{-0.5mm}
\end{figure}
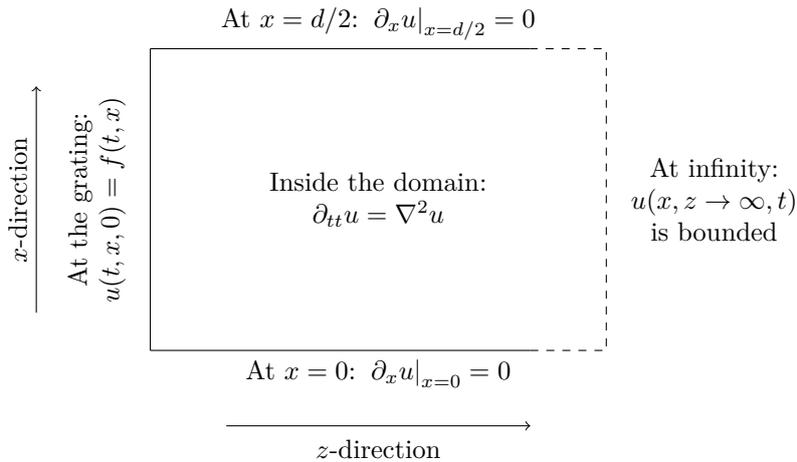

Finally, as we want to understand how the wave spreads after passing through the grating, we require zero initial conditions
\begin{equation*}
    u(t=0) = 0,\quad\text{and,}\quad \partial_t u(t=0) = 0.
\end{equation*}
By putting everything together we have the following well-posed problem:
\begin{empheq}[box=\fbox]{align}
\partial_{tt} u(t,x,z) &= \partial_{xx} u(t,x,z) + \partial_{zz} u(t,x,z) ,\quad t > 0, \, c\in \mathbb{R},\, z > 0,\nonumber \\
s.t.\hspace{2cm}&\nonumber \\
    u(t = 0, x, z) &=  \partial_t u(t=0, x, z) = 0,\label{pde_statement}\\
\partial_x u(t,x = 0, z) &=  \partial_x u(t,x = d/2, z) = 0,\nonumber \\
    u(t,x, z = 0) &= g(x) h(t)\, \theta(t),\quad u(t,x,z \to \infty) \text{ is bounded}.\nonumber 
\end{empheq}
This is summarised in Figure \ref{problem_diagram}.

There are two things to comment about this model: First, we have \emph{not} required $u(z=0) = \partial_z u(z=0)=0$ at the screen, which are the mathematically inconsistent Kirchhoff conditions. Similarly, we are not discarding functions of $t+z$, that is incoming waves, because of their unphysical nature. Secondly, and contrary to what is usually expected, our scalar field $u \equiv A^y $ is not one of the physical electromagnetic fields, but rather its potential. So $\abs{u}^2$ is not the intensity, although both are closely related.

\section{Solution of the model: The transient solution}\label{section:transient_solution}

In this section we solve the well-posed problem \eqref{pde_statement}. We start by using the Fourier decomposition of $u$ in $x$
\begin{equation}\label{eq:fourier_u}
    u(t,x,z) = \sum_{n\in\mathbb{N}} c_n(t,z) \cos{k_n x},
\end{equation}
where we have defined
\begin{equation}\label{eq:k_n}
    k_n := \frac{2\pi}{d} n.
\end{equation}
By substituting this into the wave equation, we find
\begin{equation*}
    \sum_{n\in\mathbb{N}} \partial_{tt} c_n \cos{k_n x} = \sum_{n\in\mathbb{N}} \left( - k_n^2 c_n + \partial_{zz} c_n \right) \cos{k_n x},
\end{equation*}
from which we can extract that the Fourier coefficients satisfy a 1D Klein-Gordon equation,
\begin{equation*}
    \partial_{tt} c_n - \partial_{zz}c_n + k_n^2 c_n = 0.
\end{equation*}
By combining equation \eqref{eq:fourier_u} with the boundary conditions of $u$ we can recover the boundary conditions that $c_n$ must satisfy:
\begin{align*}
    c_n(z = 0) &= g_n h(t)\, \theta(t), & c_n(z \to \infty) &\text{ is bounded},\\
    c_n(t = 0) &= 0, & \partial_t c_n(t=0) &= 0,
\end{align*}
where $\hat{g}_n = \hat{g}_{-n} = \frac{2}{d}\int_0^{d/2} g(x) \cos k_n x \,dx$ is the Fourier coefficient of $g(x)$.\smallskip

To make progress, we Laplace-transform the equation in time, so we get
\begin{align*}
    \mathcal{L}_t\left\{\partial_{tt} c_n - \partial_{zz}c_n + k_n^2 c_n\right\} (s) &= s^2 C_n - s c_n(t=0,z) - \partial_t c_n(t=0,z) - \partial_{zz}C_n + k_n^2 C_n \\
    &= s^2 C_n - \partial_{zz}C_n + k_n^2 C_n,
\end{align*}
where 
\begin{equation*}
    C_n(s,z) \equiv \mathcal{L}_t\{c_n(t,z)\}(s) \equiv \int_0^\infty c_n(t,z) e^{-ts}\,dt
\end{equation*}
is the Laplace transform of $c_n$. We thus have the equation
\begin{equation*}
    \partial_{zz}C_n - (s^2 + k_n^2) C_n = 0,
\end{equation*}
with boundary conditions
\begin{align}\label{eq:BC_z}
    C_n(s,z = 0) &= g_n H(s), & C_n(s,z \to \infty) &\text{ is bounded},
\end{align}
where $H(s) = \mathcal{L}_t\{h(t)\}(s)$. The solution to this equation is of the form
\begin{equation*}
    C_n (s,z) = A_n(s) e^{- \sqrt{s^2 + k_n^2} z} + B_n(s) e^{\sqrt{s^2 + k_n^2} z},
\end{equation*}
and by imposing the boundary conditions \eqref{eq:BC_z}, which ensure the uniqueness of the solution, we find
\begin{equation*}
    C_n (s,z) = g_n H(s) e^{- \sqrt{s^2 + k_n^2} z} .
\end{equation*}
Thus, the solution to our Klein-Gordon problem is
\begin{equation*}
    c_n (t,z) = \mathcal{L}^{-1}_s \left\{ C_n(s,z) \right\}(t) = \mathcal{L}^{-1}_s \left\{ g_n H(s) e^{- \sqrt{s^2 + k_n^2} z} \right\}(t).
\end{equation*}\smallskip

To invert the Laplace transform of $C_n$ we will make use that the Laplace transform of the convolution of two functions verifies
\begin{equation*}
    \mathcal{L}_t\{f(t)*g(t)\}(s) = F(s)G(s),
\end{equation*}
where $*$ denotes the convolution operator $f(t) * g(t) = \int_0^t f(\tau) g(t-\tau)\,d\tau$, and $F$ and $G$ are the Laplace transforms of $f$ and $g$, respectively. Using the facts that 
\begin{align}\label{eq:invLT_G}
    \mathcal{L}_s^{-1} \left\{  e^{-z\sqrt{k_n^2 + s^2}} \right\} (t) &=  \delta(t-z) - \theta(t-z) k_n z \dfrac{J_1\left( k_n\sqrt{t^2-z^2} \right)}{\sqrt{t^2-z^2}},\\
    \mathcal{L}_s^{-1} \left\{  \hat{g}_n H(s) \right\} (t) &=  \hat{g}_n h(t) \,\theta(t),\nonumber
\end{align}
where we have used equation \eqref{eq:appendix_result} and $\delta(t)$ is Dirac's delta function, we find that
\begin{align*}
\begin{split}
    c_n(t,z) &= \hat{g}_n \int_0^t h(t-\tau) \left( \delta(\tau-z) - \theta(\tau-z)\, k_n z\, \dfrac{J_1\left( k_n\sqrt{\tau^2-z^2} \right)}{\sqrt{\tau^2-z^2}} \right) \,d\tau\\
    &= \hat{g}_n h(t-z) \theta(t-z) - k_n z  \int_0^t h(t-\tau) \theta(\tau-z) \dfrac{J_1\left( k_n\sqrt{\tau^2-z^2} \right)}{\sqrt{\tau^2-z^2}} \,d\tau\\
    &= \hat{g}_n h(t-z) \theta(t-z) - \hat{g}_n k_n z \int_z^{\max(t,z)} h(t-\tau) \dfrac{J_1\left( k_n\sqrt{\tau^2-z^2} \right)}{\sqrt{\tau^2-z^2}} \,d\tau.
\end{split}
\end{align*}
Putting everything together, we find that the Talbot effect is described by the following function:
\begin{equation}\label{eq:general_solution}
    u(t,x,z) = \sum_n \hat{g}_n \left( h(t-z) - k_n z \int_z^{t} \dfrac{J_1\left( k_n\sqrt{\tau^2-z^2} \right)}{\sqrt{\tau^2-z^2}} h(t-\tau) \,d\tau \right) \theta(t-z) \cos{k_n x},
\end{equation}
which in the monochromatic case, $h(t) = \sin \omega t$, becomes
\begin{empheq}[box=\fbox]{equation}\label{eq:transient_solution}
    u(t,x,z) = \sum_n \hat{g}_n\left( \sin\omega (t-z) - k_n z \int_z^{t} \dfrac{J_1\left( k_n\sqrt{\tau^2-z^2} \right)}{\sqrt{\tau^2-z^2}} \sin\omega(t-\tau) \,d\tau \right) \theta(t-z) \cos{k_n x}.
\end{empheq}
Remember that this is only defined for $t\geq 0$, $z\geq 0$. We will call this the \emph{full solution} or the \emph{transient solution} to the problem.\smallskip

Note that at no point we have reduced our problem to a Helmholtz equation, which has saved us from needing Sommerfeld's radiation condition to choose a solution, as is the case in earlier treatments of this problem; see \cite{Eceizabarrena2017} for an explicit use in the case of the Schrödinger equation.\smallskip

\section{Recovering the stationary solution}\label{section:stationary_solution}

We can see that $u$ --the solution that we obtained in \eqref{eq:transient_solution}-- is a more complicated expression than $U$ --the one in the literature \eqref{eq:stationary_solution_cited}--, and it is not clear whether both agree. In this section we will show that $u$ evolves into $U$, in the sense that the leading asymptotic behaviour of $u$ as $t\gg z$ is precisely $ \text{Im} \left[U(x,z) e^{i\omega t} \right]$. The importance of this result is that it allows us to mathematically justify the use of equation \eqref{eq:stationary_solution_cited}, which is now rigorously derived from the wave equation. More formally, we prove the following statement.
\begin{theo}
    Given the solution to problem \eqref{pde_statement},
    \begin{equation*}
        u(t,x,z) = \sum_n \hat{g}_n\left( \sin\omega (t-z) - k_n z \int_z^{t} \dfrac{J_1\left( k_n\sqrt{\tau^2-z^2} \right)}{\sqrt{\tau^2-z^2}} \sin\omega(t-\tau) \,d\tau \right) \theta(t-z) \cos{k_n x}, \quad k_n = \frac{2\pi }{d}n,
    \end{equation*}
    and the solution to the Helmholtz equation
    \begin{equation}\label{eq:U_definition}
        U(x,z) = \sum_n \hat{g}_n \cos{k_n x}\times \begin{cases}
           \exp \left(- iz\sqrt{\omega^2 - k_n^2} \right), & \quad 0< k_n \leq\omega \\ 
           \exp \left(-z\sqrt{k_n^2 - \omega^2}\right), & \quad \omega< k_n <\infty,
        \end{cases}
    \end{equation}
    we have that the asymptotic behaviour of $u$ is precisely $\text{Im} \left[U(x,z) e^{i\omega t} \right]$. More precisely,
    \begin{equation}\label{eq:u_asymptotic_limit}
        \lim_{t\to 0} u(t,x,z) - \text{Im} \left[U(x,z) e^{i\omega t} \right] = 0.
    \end{equation}
\end{theo}\smallskip

\emph{Proof:} In order to show this, we will rewrite the integral as
\begin{equation*}
    \int_z^{t} \sin\omega(t-\tau) \dfrac{J_1\left( k_n\sqrt{\tau^2-z^2} \right)}{\sqrt{\tau^2-z^2}} \,d\tau = \left[ \int_z^{\infty} - \int_t^{\infty} \right] \sin\omega(t-\tau) \dfrac{J_1\left( k_n\sqrt{\tau^2-z^2} \right)}{\sqrt{\tau^2-z^2}} \,d\tau,
\end{equation*}
and study the leading integral. By making the change of variable $r^2 = \tau^2 - z^2$ we get
\begin{equation*}
    \int_z^{\infty} \sin\omega(t-\tau) \dfrac{J_1\left( k_n\sqrt{\tau^2-z^2} \right)}{\sqrt{\tau^2-z^2}} \,d\tau = \int_0^{\infty} \dfrac{\sin\omega(t-\sqrt{r^2+z^2})}{\sqrt{r^2+z^2}} J_1\left( k_n r \right) \,d\tau,
\end{equation*}
which can be rewritten into a Hankel integral:
\begin{align*}
\begin{split}
    I(t,z;k) := \int_0^{\infty} \dfrac{\sin\omega(t-\sqrt{r^2+z^2})}{\sqrt{r^2+z^2}} J_1\left( k_n r \right) \,d\tau &= \int_0^{\infty} \dfrac{\sin\omega(t-\sqrt{r^2+z^2})}{\sqrt{k_n r} \sqrt{r^2+z^2}} \sqrt{k_n r} J_1\left( k_n r \right) \,d\tau \\
    &\equiv \mathcal{H}_1 \left(\dfrac{\sin\omega(t-\sqrt{r^2+z^2})}{\sqrt{k_n r} \sqrt{r^2+z^2}};k_n \right) .
\end{split}
\end{align*}
Now we use the fact that 
\begin{equation*}
    \sin\omega \left(t-\sqrt{r^2+z^2}\right) = \sin\omega t\, \cos(\omega\sqrt{r^2+z^2}) - \cos\omega t\, \sin(\omega\sqrt{r^2+z^2}),
\end{equation*}
and the following tabulated Hankel transforms from \cite{bateman1954tables}\\
\begin{center}
\centering
\begin{tabular}{cc}
\toprule
$f(r)$ & $\mathcal{H}_1 \left(f(r);k \right)  \equiv \int_0^{\infty} f(r) \sqrt{k_n r} J_1\left( k_n r \right) \,dr$ \\\midrule
$\dfrac{\sin\omega \sqrt{r^2+z^2}}{\sqrt{r} \sqrt{r^2+z^2}}$ & $\begin{cases}
  \frac{1}{z\sqrt{k}} \left[ \sin\omega z - \sin z\sqrt{\omega^2-k^2} \right],  & 0 < k < \omega \\
  \frac{\sin\omega z}{z\sqrt{k}}, & \omega < k < \infty
\end{cases}$ \\
$\dfrac{\cos\omega \sqrt{r^2+z^2}}{\sqrt{r} \sqrt{r^2+z^2}}$ & $\begin{cases}
  \frac{1}{z\sqrt{k}} \left[ \cos\omega z - \cos z\sqrt{\omega^2-k^2} \right],  & 0<k<\omega \\
  \frac{1}{z\sqrt{k}}\left[ \cos\omega z - e^{-z\sqrt{k^2 - \omega^2}} \right], &  \omega<k<\infty
\end{cases}$\\\bottomrule
\end{tabular}    
\end{center}
So we find that
\begin{align*}
\begin{split}
    I(t,z;k_n) &= \sin\omega t\,\mathcal{H}_1 \left(\dfrac{\cos\omega(\sqrt{r^2+z^2})}{\sqrt{r} \sqrt{r^2+z^2}};k_n \right) - \cos\omega t\,\mathcal{H}_1 \left(\dfrac{\sin\omega(\sqrt{r^2+z^2})}{\sqrt{r} \sqrt{r^2+z^2}};k_n \right)\\
    &= \frac{1}{k_n z} \begin{cases}
           \sin\omega t\,\left[ \cos\omega z - \cos z\sqrt{\omega^2 - k_n^2} \right] - \cos\omega t\,\left[ \sin\omega z - \sin z\sqrt{\omega^2 - k_n^2} \right], & 0 < k_n \leq\omega \\
          \sin\omega t\, \left[\cos\omega z- e^{-z\sqrt{k_n^2 - \omega^2}} \right] - \cos\omega t \, \sin\omega z, &  \omega< k_n <\infty
        \end{cases}\\
    &= \frac{1}{k_n z} \begin{cases}
       \sin\omega (t - z) - \sin \left(\omega t - z\sqrt{\omega^2 - k_n^2} \right), & 0< k_n \leq\omega \\ \sin\omega (t - z) - 
        \sin\omega t\, e^{-z\sqrt{k_n^2 - \omega^2}}, & \omega< k_n <\infty.
    \end{cases}\\
\end{split}
\end{align*}
Thus we find that the asymptotic behaviour of our solution as $t\to\infty$ is 
\begin{align*}
\begin{split}
    u (t,x,z) &= \sum_n \hat{g}_n\left( \sin\omega (t-z) - k_n z \left[\int_z^{\infty} - \int_t^{\infty} \right]\dfrac{J_1\left( k_n\sqrt{\tau^2-z^2} \right)}{\sqrt{\tau^2-z^2}} \sin\omega(t-\tau) \,d\tau \right) \theta(t-z) \cos{k_n x}\\
    &= \sum_n \hat{g}_n \cos{k_n x}\times \begin{cases}
       \sin \left(\omega t - z\sqrt{\omega^2 - k_n^2} \right) + E_n(t,z), &  0< k_n \leq\omega \\ 
       \sin\omega t\, e^{-z\sqrt{k_n^2 - \omega^2}} + E_n(t,z), &  \omega< k_n <\infty,
    \end{cases}
\end{split}
\end{align*}
where we have defined the error
$$
E_n(t,z) := \theta(t-z) \,k_n z \int_t^{\infty} \dfrac{J_1\left( k_n\sqrt{\tau^2-z^2} \right)}{\sqrt{\tau^2-z^2}} \sin\omega(t-\tau) \,d\tau .
$$
We can show that the error is of order $E_n = O\left(\sqrt{\frac{z}{t}}\right)$, which we do in Appendix \ref{appendix:bounding_E}. Thus,

\begin{empheq}[box=\fbox]{equation}\label{eq:stationary_solution}
    u (t,x,z) = \sum_n \hat{g}_n \cos{k_n x}\times \begin{cases}
       \sin \left(\omega t - z\sqrt{\omega^2 - k_n^2} \right) + O\left(\sqrt{\dfrac{z}{t}}\right), &  0< k_n \leq\omega \\ 
       \sin\omega t\, e^{-z\sqrt{k_n^2 - \omega^2}} + O\left(\sqrt{\dfrac{z}{t}}\right), &  \omega< k_n <\infty,
    \end{cases}
\end{empheq}
and thus
\begin{equation*}
    u (t,x,z) - \text{Im} \left[U(x,z) e^{i\omega t} \right] = \sum_n \hat{g}_n \cos{k_n x}\times  O\left(\sqrt{\dfrac{z}{t}}\right),
\end{equation*}
from which \eqref{eq:u_asymptotic_limit} follows.\hfill $\square$

\section{Talbot carpets}\label{section:graphing_carpets}

\subsection{Graphing the solution}
Now that we have obtained a function \eqref{eq:transient_solution} that describes the Talbot effect we want to check whether this agrees with Lord Talbot's findings and study any further properties. The first thing that we do is to plot the intensity of the function, that is the square of the function. For this, we first need to pick the shape of the grating. We will chose a Ronchi grating, that is, a perfect slit of width $l$ centered at $x=0$ which lets through light from a source with amplitude $A$. Thus
\begin{equation*}\label{eq:rect_coeffs}
    g(x) := A \dfrac{d}{l}\,\chi\left(\frac{x}{w}\right) = \begin{cases}
        A\dfrac{d}{l}, & x < \frac{x}{w} \\
        0, & x > \frac{x}{w},
    \end{cases} \quad\text{and}\quad \hat{g}_n = A \dfrac{d}{l} \hat{\chi}_n = \begin{cases}
        A, & n = 0 \\
        A \dfrac{d}{l} \dfrac{\sin{\frac{n\pi l}{d}}}{n \pi}, & n \neq 0.
    \end{cases}
\end{equation*}
Note that we have introduced a $\frac{d}{l}$ coefficient so that no energy is absorbed by the grating, and to make it converge to the Dirac comb as $l/d \to 0$. We also define the \emph{wavelength} $\lambda$ and the \emph{Talbot distance} $z_T$ as
\begin{align}
    \lambda :=& \dfrac{2\pi}{\omega} & z_T := \dfrac{2 d^2}{\lambda}= \dfrac{d^2 \omega}{\pi}.
\end{align}\smallskip

Now, the sum must be truncated. For this, it is important to note that all the terms with $k_n < \omega$ will asymptotically lead to an oscillating term, which should be accounted for, as they do not decay while propagating away from the grating. In contrast, the terms with $k_n > \omega$ will decay exponentially as they propagate in the $z$-direction, so we can neglect most of them. We cannot neglect all as they may still be large enough to matter in the proximity of $z=0$.Because of this, we have chosen to sum up to $5\,d/\lambda$ terms to ensure that the factor $1/n$ becomes small enough. This also means there will be four times as many evanescent waves as oscillating ones.\smallskip

\begin{figure}[H]
\centering
    \begin{subfigure}{.45\textwidth}
      \centering
      \includegraphics[width=\linewidth]{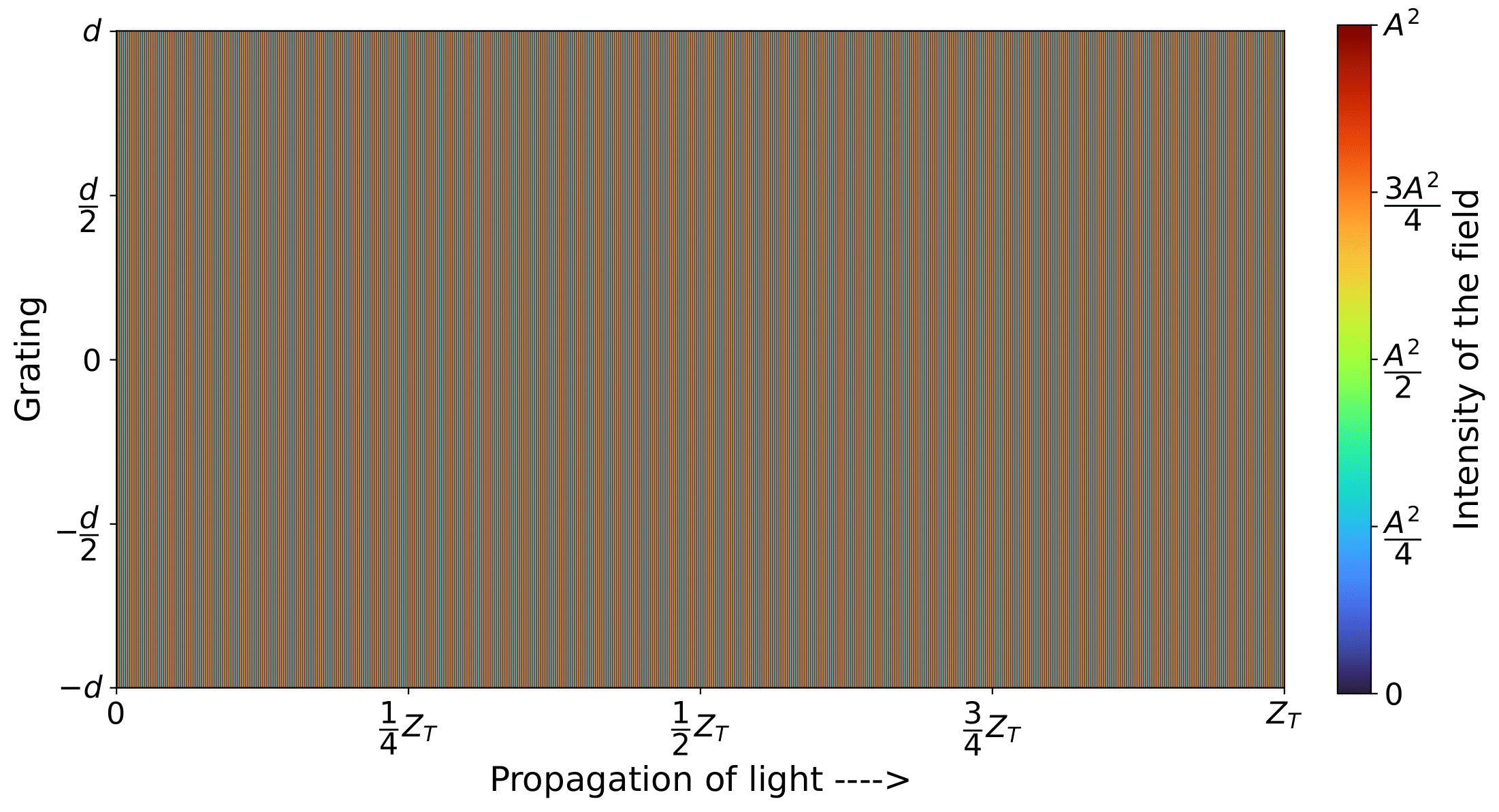}
      \caption{Carpet with $\frac{d}{\lambda}=5$ and $\frac{l}{\lambda}=5$.}
      \label{fig:d_lambda=5_w_lambda=5_carpet}
    \end{subfigure}
    \begin{subfigure}{.45\textwidth}
      \centering
      \includegraphics[width=\linewidth]{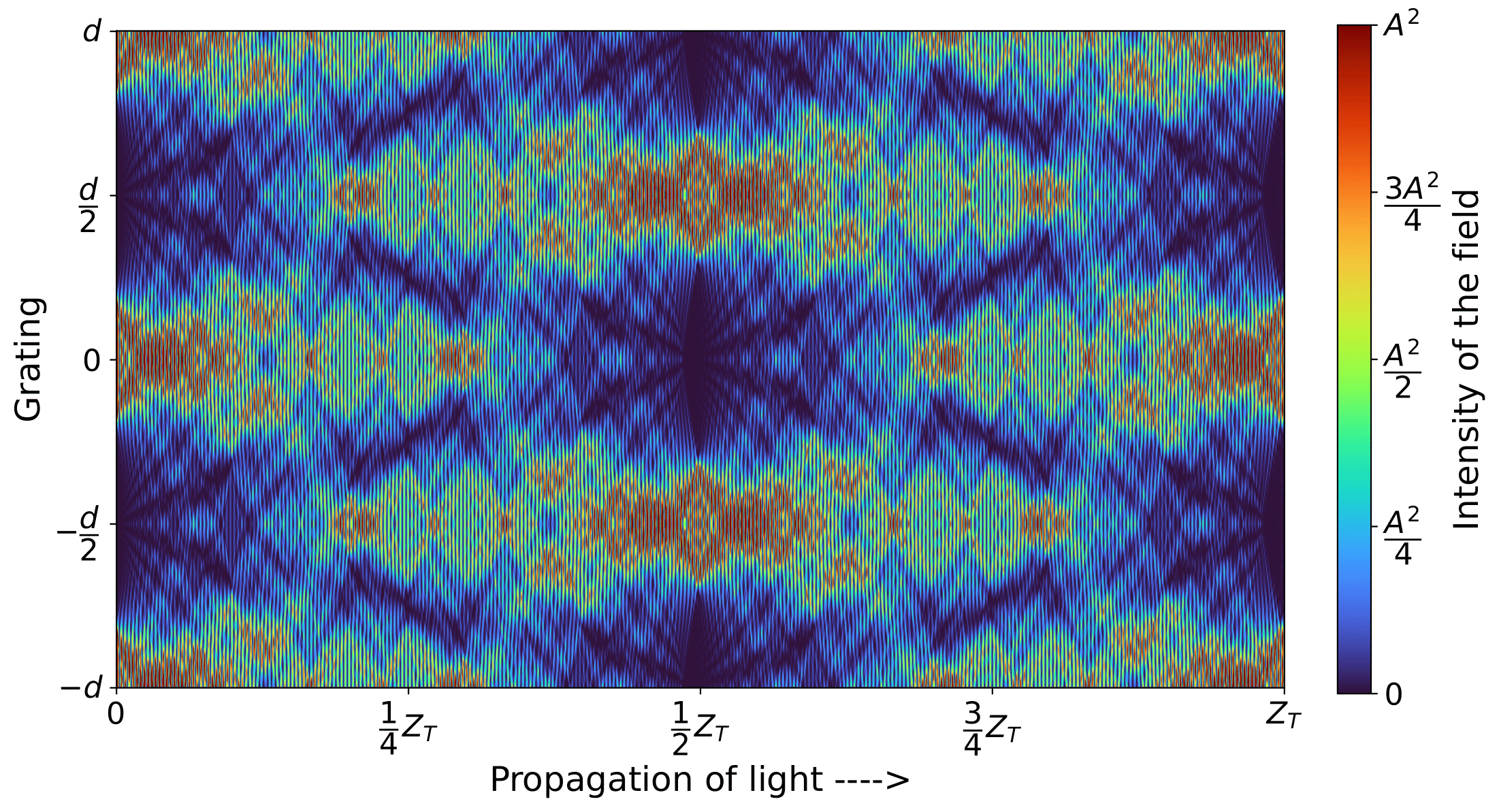}
      \caption{Carpet with $\frac{d}{\lambda}=5$ and $\frac{l}{\lambda}=2$.}
      \label{fig:d_lambda=5_w_lambda=2_carpet}
    \end{subfigure}
    
    \begin{subfigure}{.45\textwidth}
      \centering
      \includegraphics[width=\linewidth]{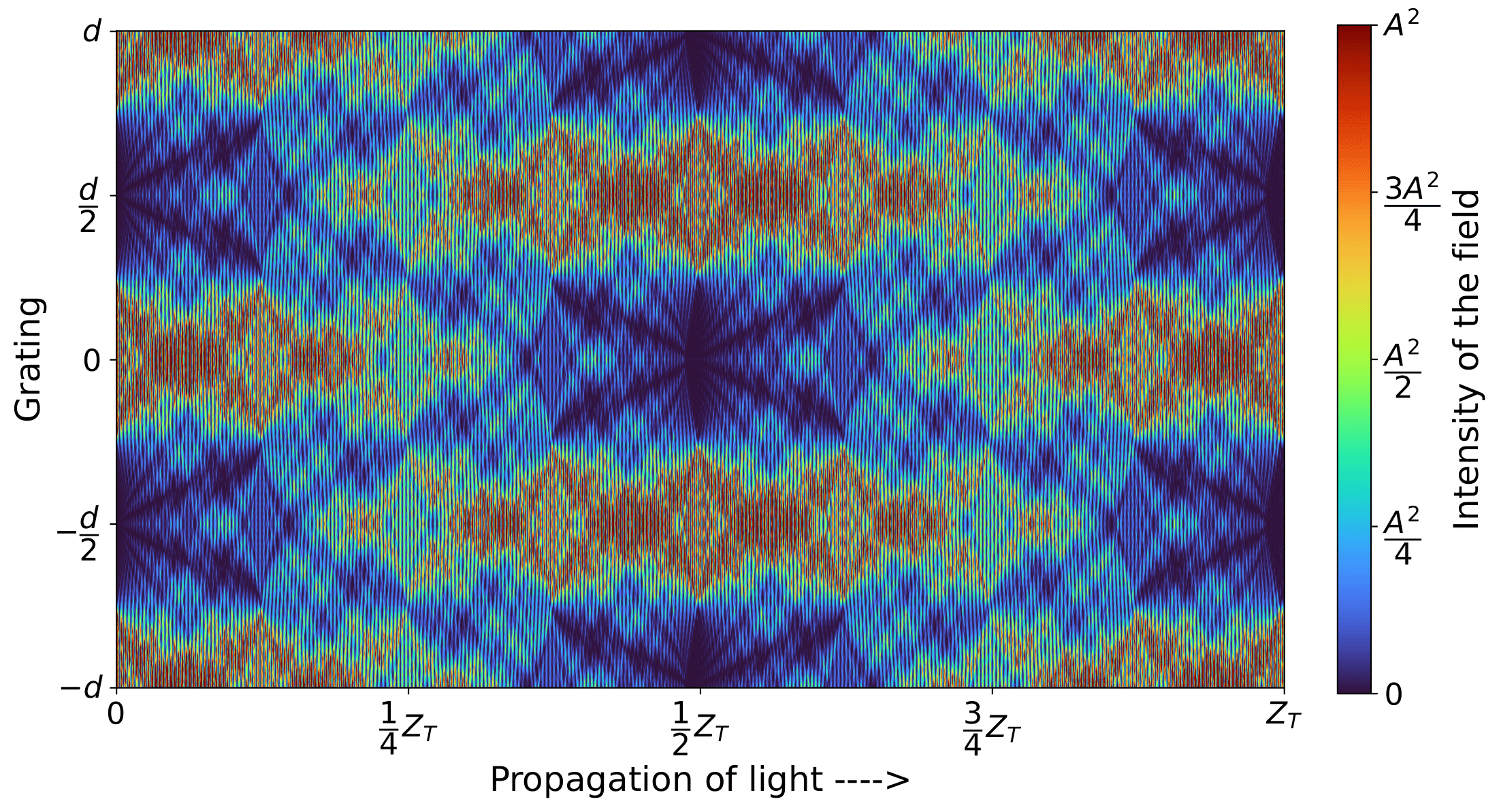}
      \caption{Carpet with $\frac{d}{\lambda}=10$ and $\frac{l}{\lambda}=5$.}
      \label{fig:d_lambda=10_w_lambda=5_carpet}
    \end{subfigure}
    \begin{subfigure}{.45\textwidth}
      \centering
      \includegraphics[width=\linewidth]{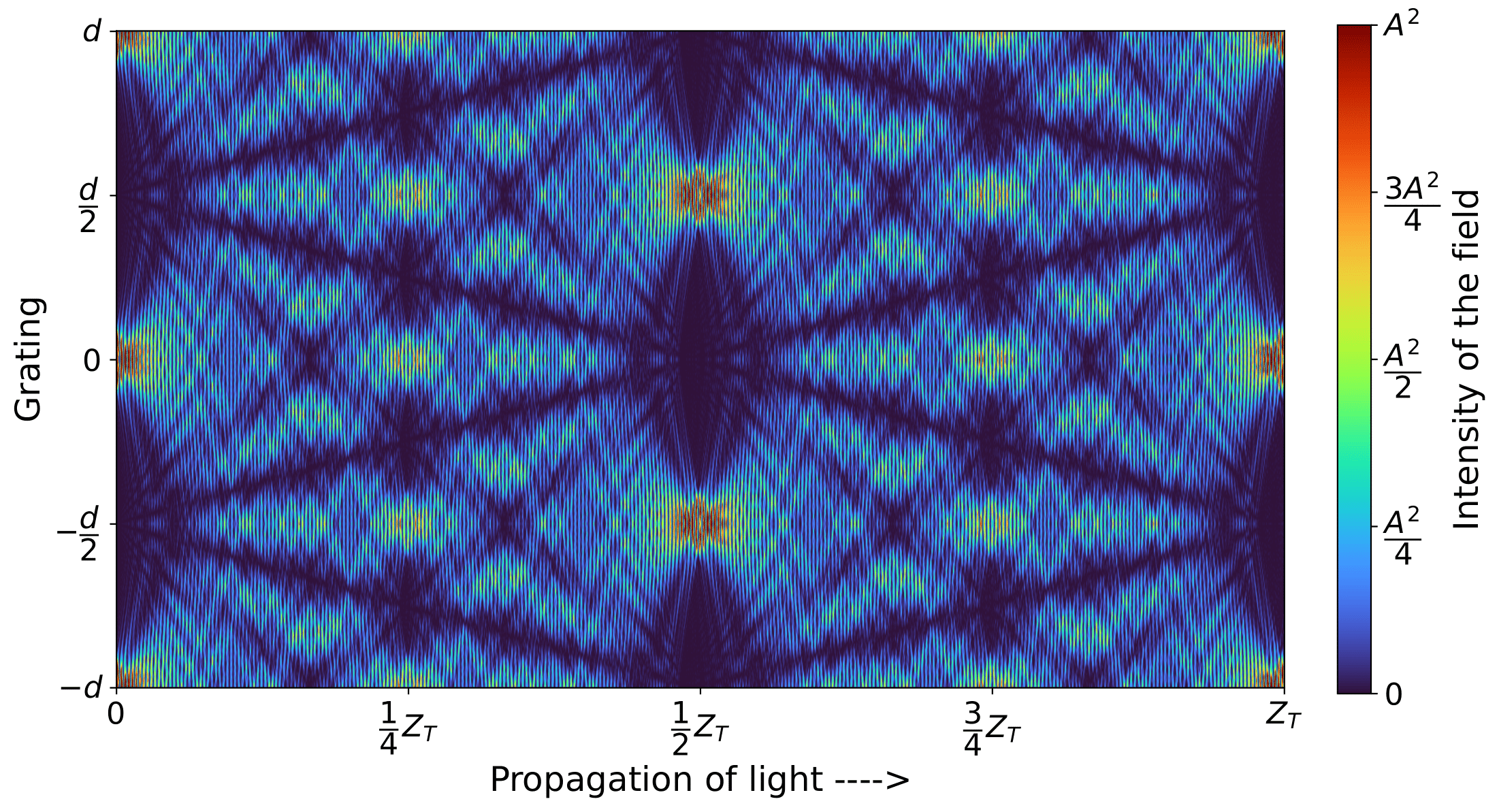}
      \caption{Carpet with $\frac{d}{\lambda}=10$ and $\frac{l}{\lambda}=2$.}
      \label{fig:d_lambda=10_w_lambda=2_carpet}
    \end{subfigure}
    
    \begin{subfigure}{.45\textwidth}
      \centering
      \includegraphics[width=\linewidth]{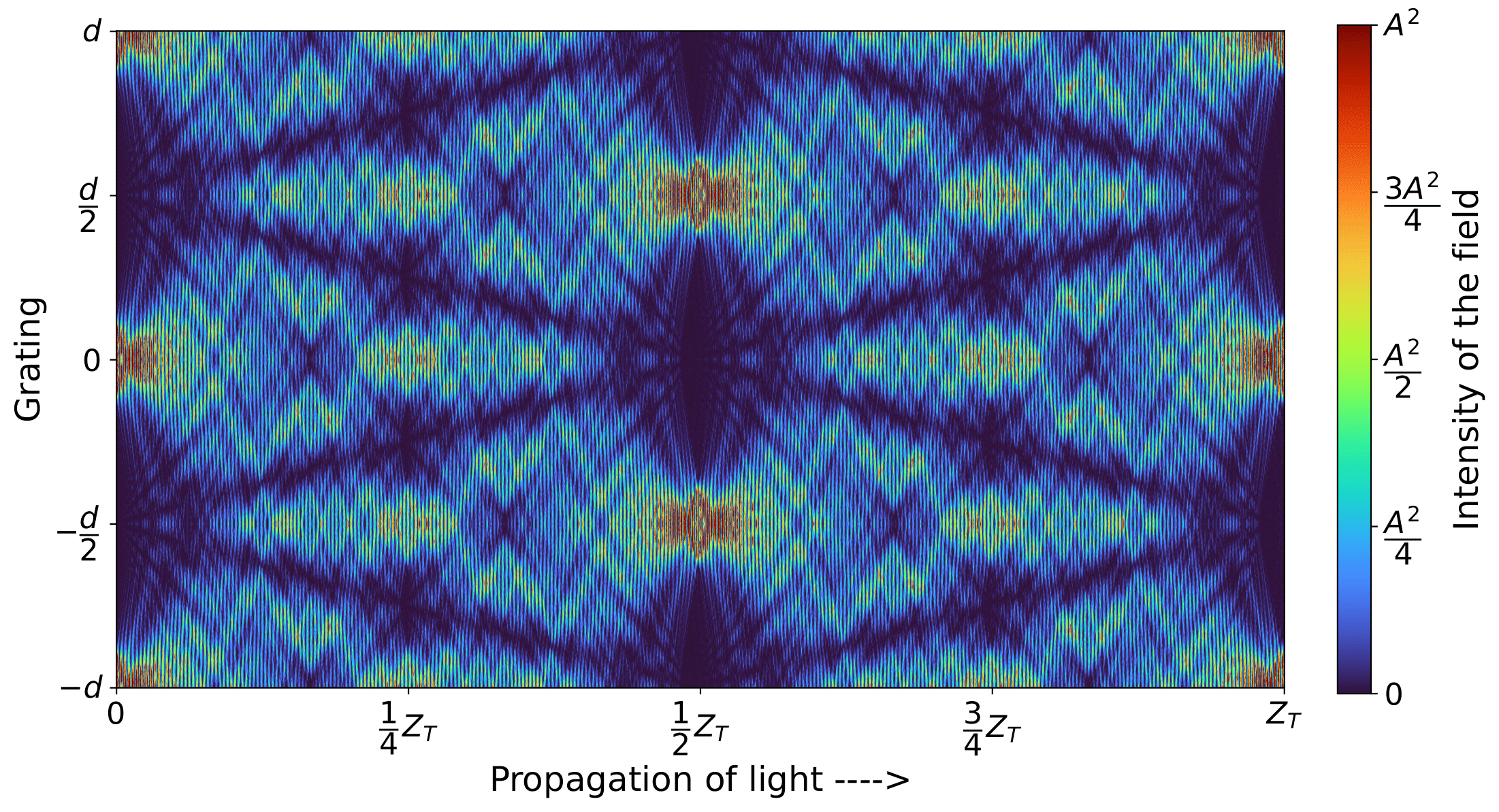}
      \caption{Carpet with $\frac{d}{\lambda}=20$ and $\frac{l}{\lambda}=5$.}
      \label{fig:d_lambda=20_w_lambda=5_carpet}
    \end{subfigure}
    \begin{subfigure}{.45\textwidth}
      \centering
      \includegraphics[width=\linewidth]{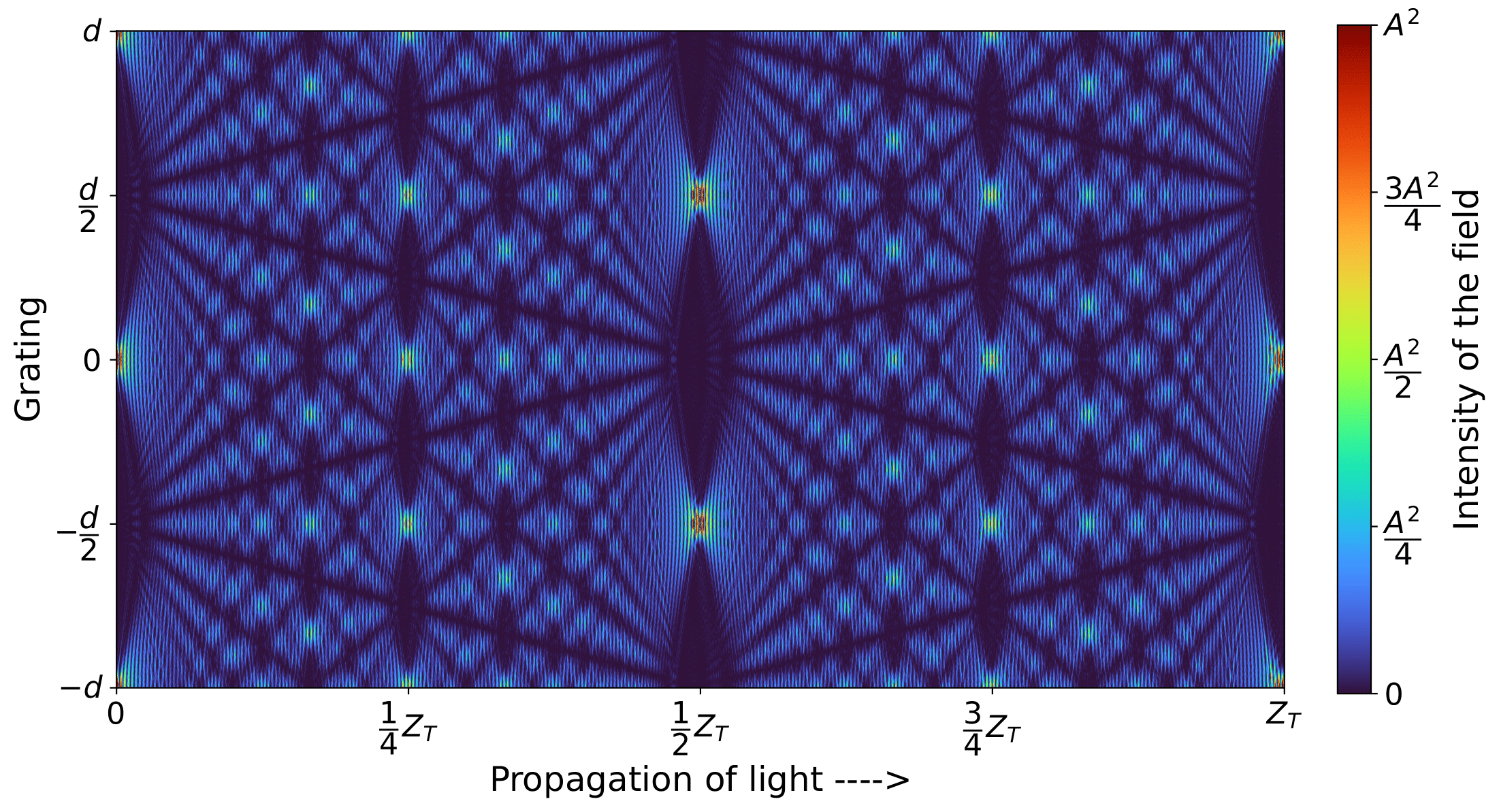}
      \caption{Carpet with $\frac{d}{\lambda}=20$ and $\frac{l}{\lambda}=2$.}
      \label{fig:d_lambda=20_w_lambda=2_carpet}
    \end{subfigure}

    \begin{subfigure}{.45\textwidth}
      \centering
      \includegraphics[width=\linewidth]{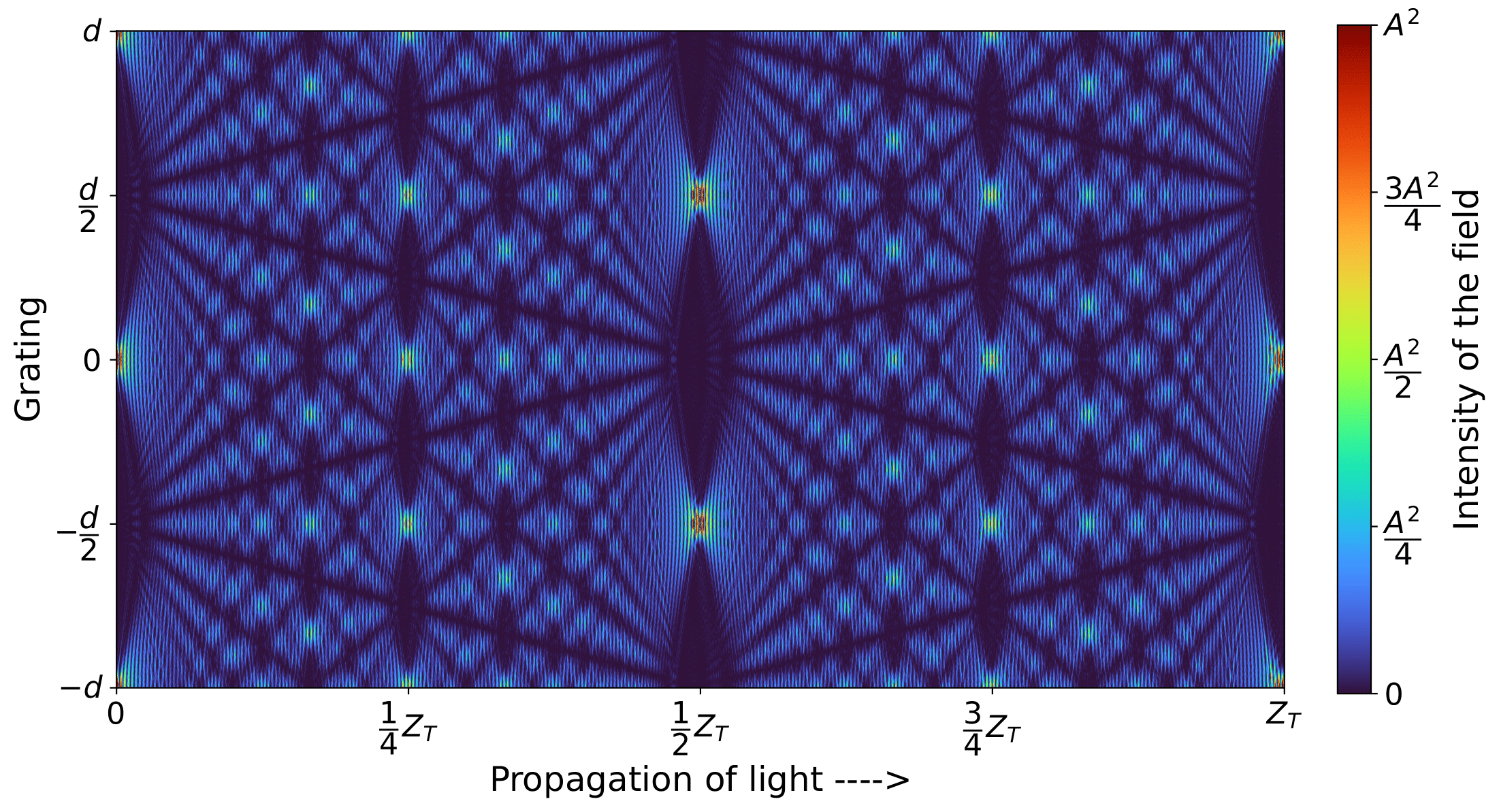}
      \caption{Carpet with $\frac{d}{\lambda}=50$ and $\frac{l}{\lambda}=5$.}
      \label{fig:d_lambda=50_w_lambda=5_carpet}
    \end{subfigure}
    \begin{subfigure}{.45\textwidth}
      \centering
      \includegraphics[width=\linewidth]{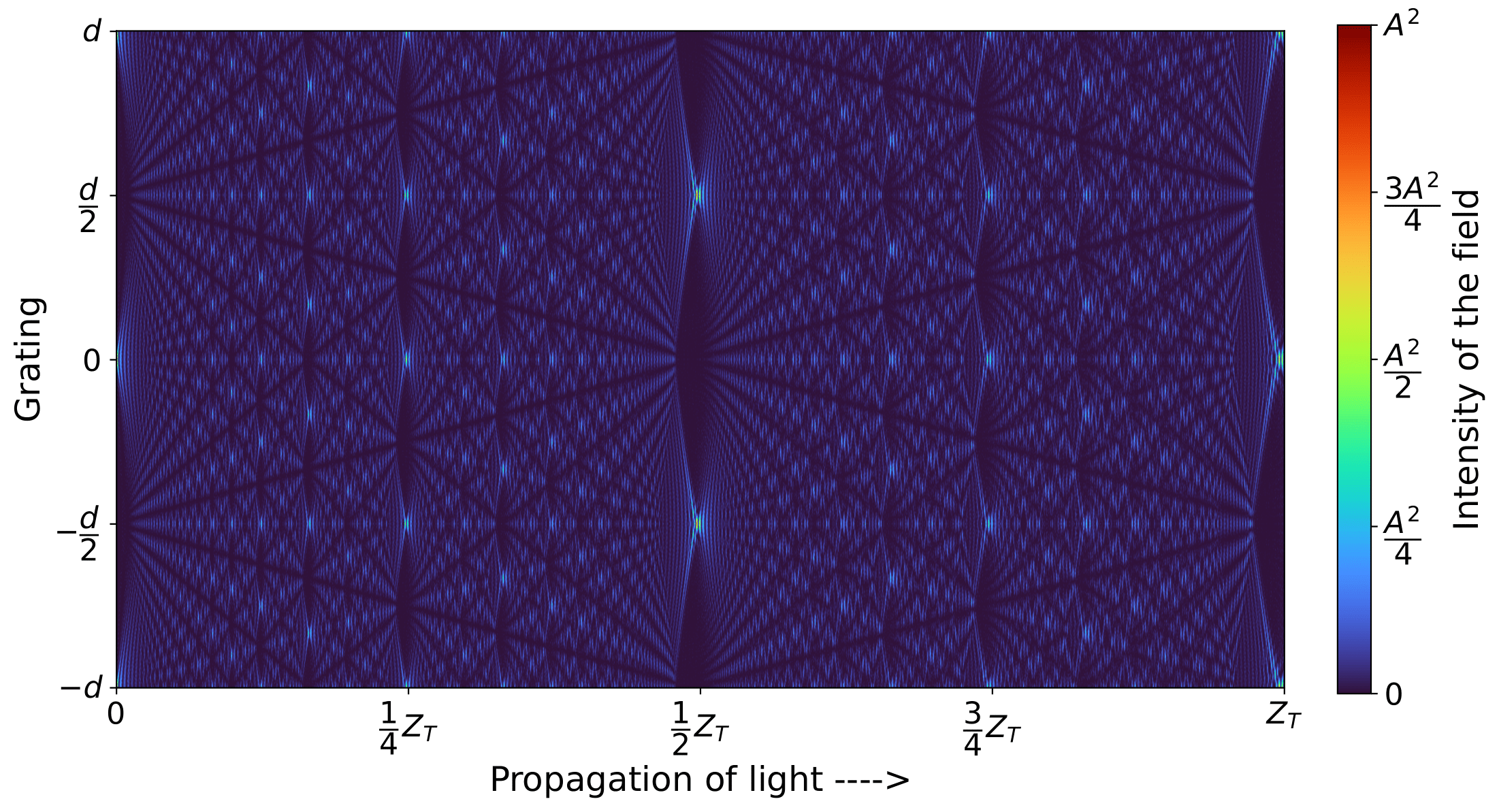}
      \caption{Carpet with $\frac{d}{\lambda}=50$ and $\frac{l}{\lambda}=2$.}
      \label{fig:d_lambda=50_w_lambda=2_carpet}
    \end{subfigure}

    \begin{subfigure}{.45\textwidth}
      \centering
      \includegraphics[width=\linewidth]{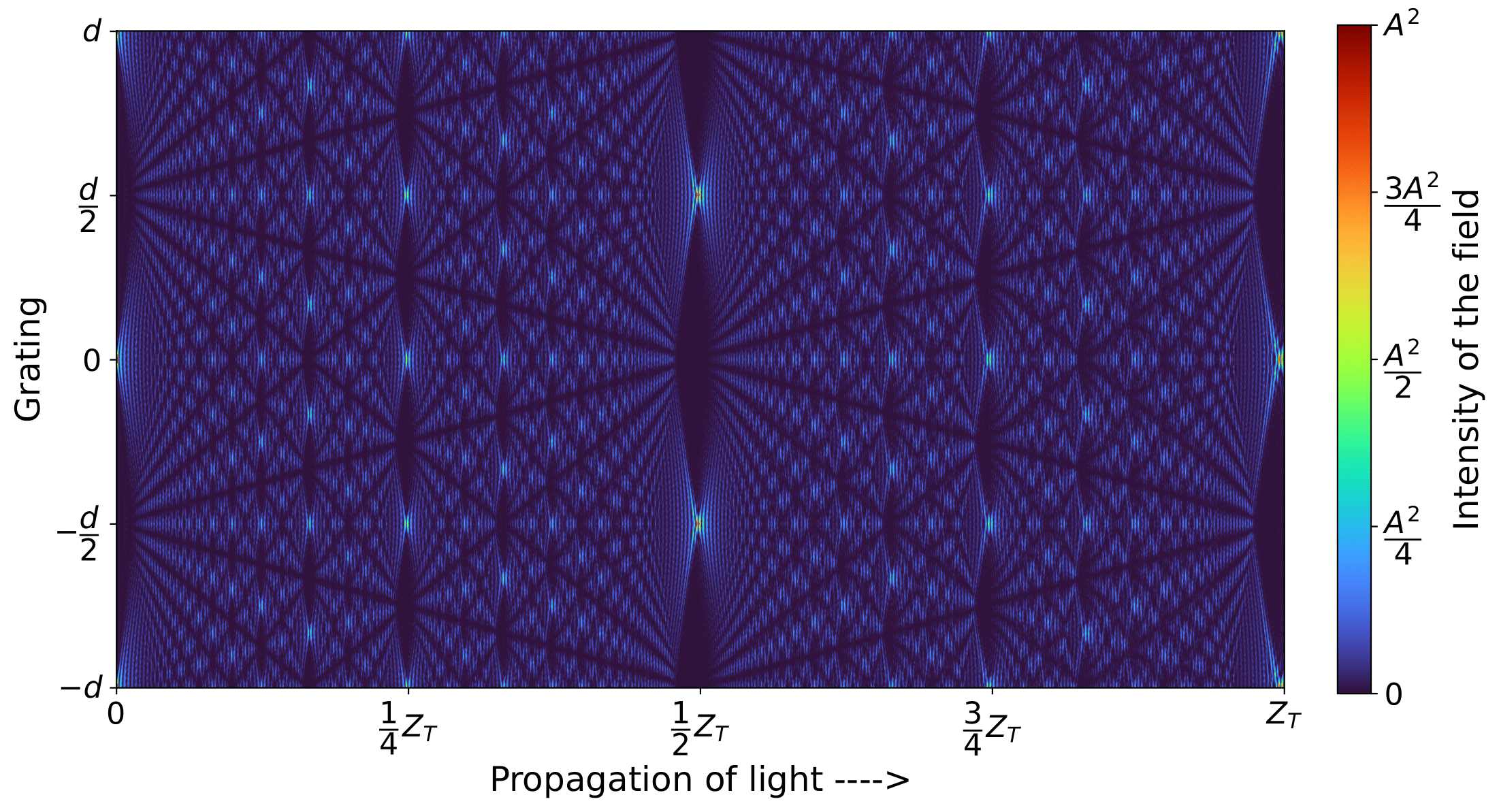}
      \caption{Carpet with $\frac{d}{\lambda}=100$ and $\frac{l}{\lambda}=5$.}
      \label{fig:d_lambda=100_w_lambda=5_carpet}
    \end{subfigure}
    \begin{subfigure}{.45\textwidth}
      \centering
      \includegraphics[width=\linewidth]{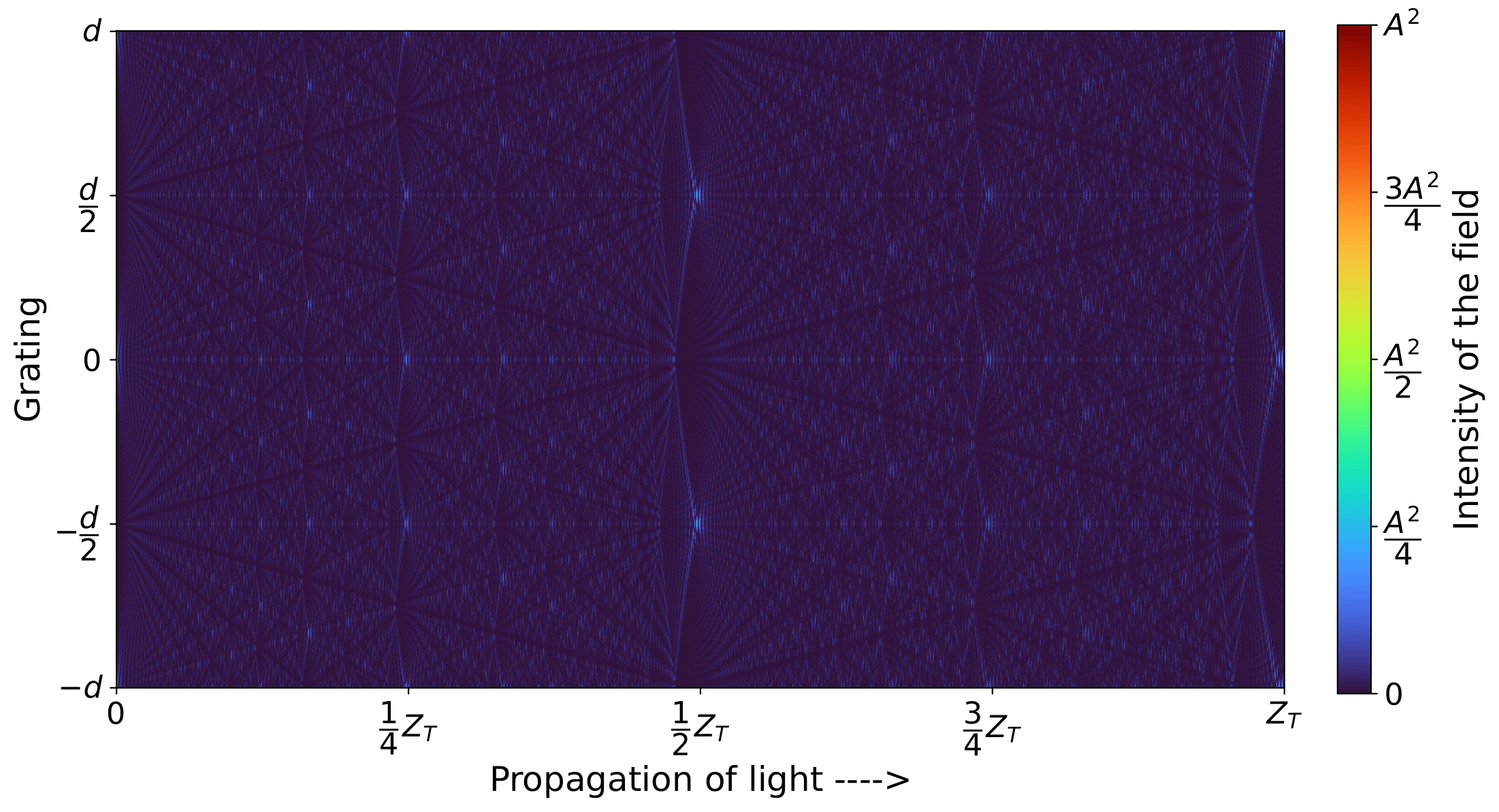}
      \caption{Carpet with $\frac{d}{\lambda}=100$ and $\frac{l}{\lambda}=2$.}
      \label{fig:d_lambda=100_w_lambda=2_carpet}
    \end{subfigure}
    
    \caption{Multiple Talbot carpets obtained by squaring and truncating equation \eqref{eq:transient_solution} at the $N = 5\,\frac{d}{\lambda}$-th term, as detailed in Section \ref{section:graphing_carpets}. A warmer colour means higher intensity. Figure \ref{fig:d_lambda=5_w_lambda=5_carpet} is just a propagating wave in free space.}
    \label{fig:multiple_carpets_transient}\vspace{-4mm}
\end{figure}

We have implemented a Python code to perform these calculations, which can be found in \cite{gabriel_ybarra_marcaida_pytalbot_2025}. With it we may create some videos that showcase the propagation of the waves, and how the Talbot effect rises. A collection of these videos for different choices of $\omega$ and $l$ can be found in \cite{YbarraVideos2025}, and captions of the final states after a time $2 z_T$ has elapsed can be found in Figure \ref{fig:multiple_carpets_transient}.
\smallskip

First of all, it can be seen that the carpets in Figures \ref{fig:multiple_carpets_transient} clearly show the Talbot effect. Indeed, in the pattern at the gratings approaches the one at the Talbot distance $z_T$. Also, fractional subimages can be seen in all these figures, being more sharp in the smaller wavelength plots. This shows that the proposed model matches the expected behaviour from the experimental findings. Moreover, thanks to the results from Section \ref{section:stationary_solution} we know that these carpets will match those obtained with the stationary solution \eqref{eq:stationary_solution}. A more detailed discussion on the agreement between the carpets resulting from \eqref{eq:stationary_solution} and experiments can be found in \cite{Berry1996, Case2009}. \smallskip

It might seem surprising that the number of the Talbot subimages increases when the parameter $\lambda/d$ approaches zero, as the effect of interferences can usually be disregarded in the short wavelength limit, as we fall back into geometric optics. But this is clearly not the case here, as the Talbot effect precisely results from the interference of waves. This apparent paradox is solved by noticing that the Talbot distance increases with $d/\lambda$, so we are dealing with interference over very long distances.\smallskip

Let us now focus on which further properties we can extract from these images. First, it can be seen how changing the wavelength and the grating's width impacts the resultant carpet. Indeed, as we decrease the ratios $d/\lambda$ and $d/w$, the carpet gets a richer structure in the sense that more rational subimages are seen. Indeed, in Figure \ref{fig:d_lambda=5_w_lambda=2_carpet} we are barely able to identify the subimages up to $z=z_T/8$ while in Figure \ref{fig:d_lambda=20_w_lambda=2_carpet} those and many more subimages are easily seen. This is also an experimental fact, and it is shown in \cite{Case2009} that subimages up to $z=z_T/12$ have been observed using $d/\lambda \approx 376$.\smallskip

To understand why we find this, we take a look at equation \eqref{eq:stationary_solution} and how changing the parameters $\lambda/d$ and $l/d$ has an impact on it. We find that diminishing the ratio $\lambda/d$ has two consequences: we increase the number of oscillating terms and we increase the number of terms for which the \emph{paraxial approximation}, that is, $\omega \gg k_n$ is true. Instead, as we decrease the ratio $l/d$ we increase the number of terms for which $\hat{g}_n \approx1$. From here we conclude that the Talbot effect arises from the paraxial terms with $\hat{g}_n \approx 1$. The rest is only noise.

\subsection{An abstraction: The amplitude envelope \texorpdfstring{$U$}{U}}

In order to make these and further properties more apparent, it is useful to strip the stationary pattern of its time evolution. We can do this by considering the envelope of the stationary pattern. Of course, this is no other than the Helmholtz' equation solution \eqref{eq:U_definition},
\begin{equation*}
    U(x,z) = \sum_n \hat{g}_n \times e^{i k_n x} \begin{cases}
       \exp \left(- iz\sqrt{\omega^2 - k_n^2} \right), &  0< k_n \leq\omega \\ 
       \exp \left(-z\sqrt{k_n^2 - \omega^2}\right), &  \omega< k_n <\infty.
    \end{cases}
\end{equation*}
This function gives us the maximum amplitude of $u$ at every point in space and thus we can easily recover the energy density over a period in the $x$-direction, which is
\begin{equation*}
    E(z) = \dfrac{1}{d} \int_0^d \abs{U (x,z)}^2 \, dx,
\end{equation*}
and because of the orthogonality of the terms of the series the cross-terms vanish, so we find
\begin{equation*}
    E(z) = \sum_n \abs{\hat{g}_n}^2 \begin{cases}
       1, &  0< k_n \leq\omega \\ 
       \exp \left(- 2 z\sqrt{k_n^2 - \omega^2}\right), &  \omega< k_n <\infty,
    \end{cases}
\end{equation*}
where we can see that some energy is dissipated as the evanescent waves travel along the $z$-direction.\smallskip

\begin{figure}[H]
\centering
    \begin{subfigure}{.45\textwidth}
      \centering
      \includegraphics[width=\linewidth]{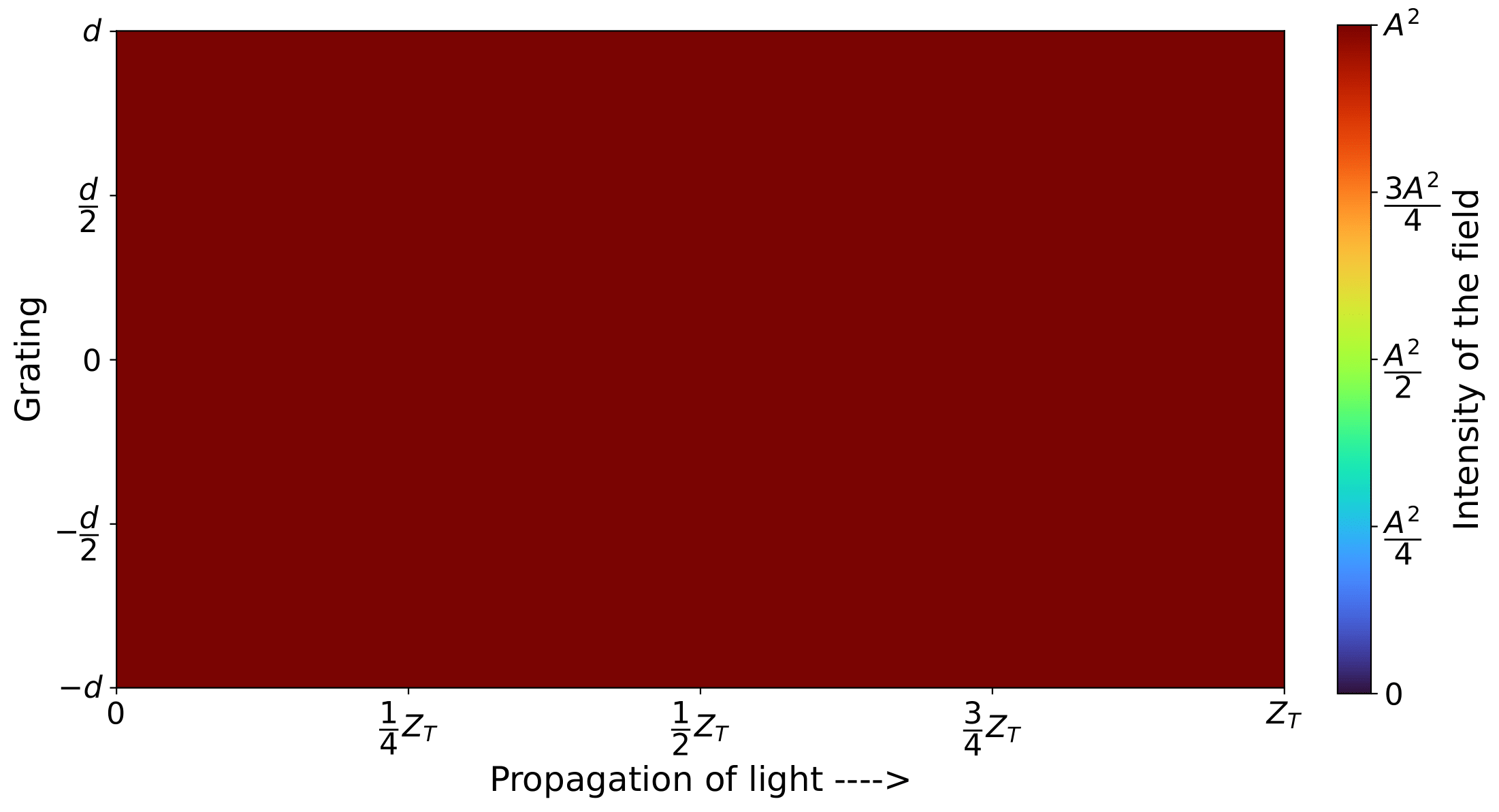}
      \caption{Carpet with $\frac{d}{\lambda}=5$ and $\frac{l}{\lambda}=5$.}
      \label{fig:d_lambda=5_w_lambda=5_energy_carpet}
    \end{subfigure}
    \begin{subfigure}{.45\textwidth}
      \centering
      \includegraphics[width=\linewidth]{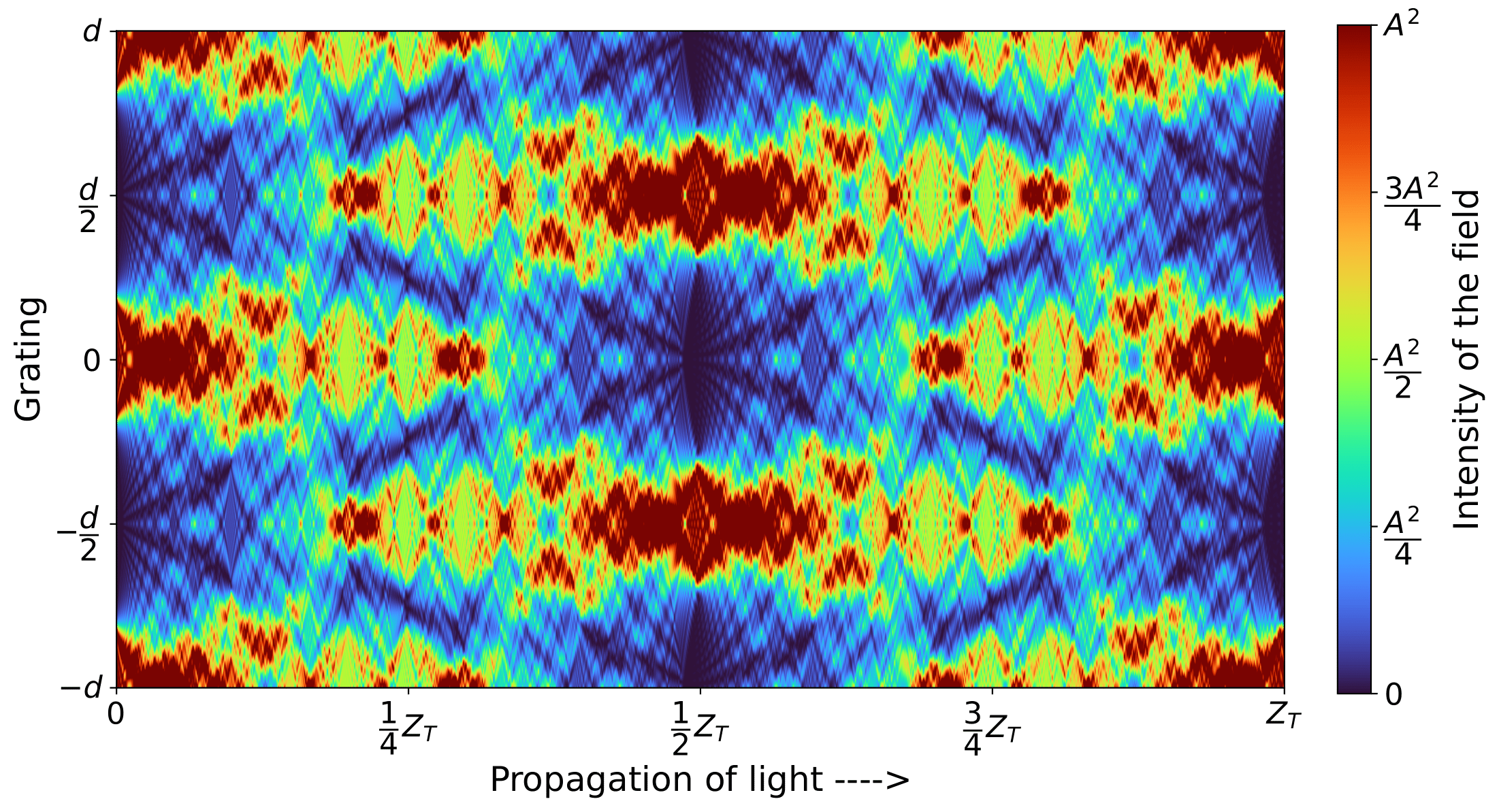}
      \caption{Carpet with $\frac{d}{\lambda}=5$ and $\frac{l}{\lambda}=2$.}
      \label{fig:d_lambda=5_w_lambda=2_energy_carpet}
    \end{subfigure}
    
    \begin{subfigure}{.45\textwidth}
      \centering
      \includegraphics[width=\linewidth]{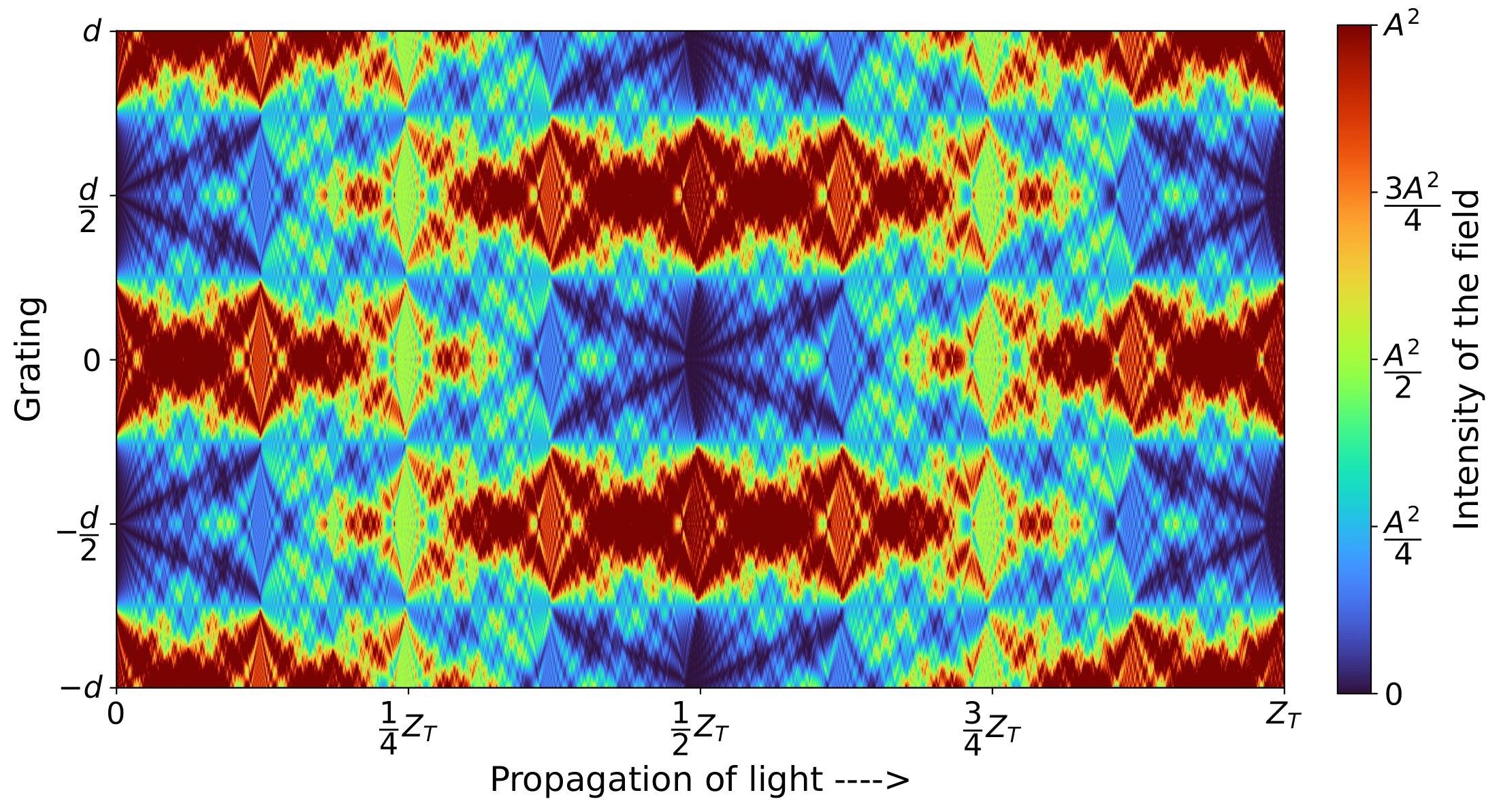}
      \caption{Carpet with $\frac{d}{\lambda}=10$ and $\frac{l}{\lambda}=5$.}
      \label{fig:d_lambda=10_w_lambda=5_energy_carpet}
    \end{subfigure}
    \begin{subfigure}{.45\textwidth}
      \centering
      \includegraphics[width=\linewidth]{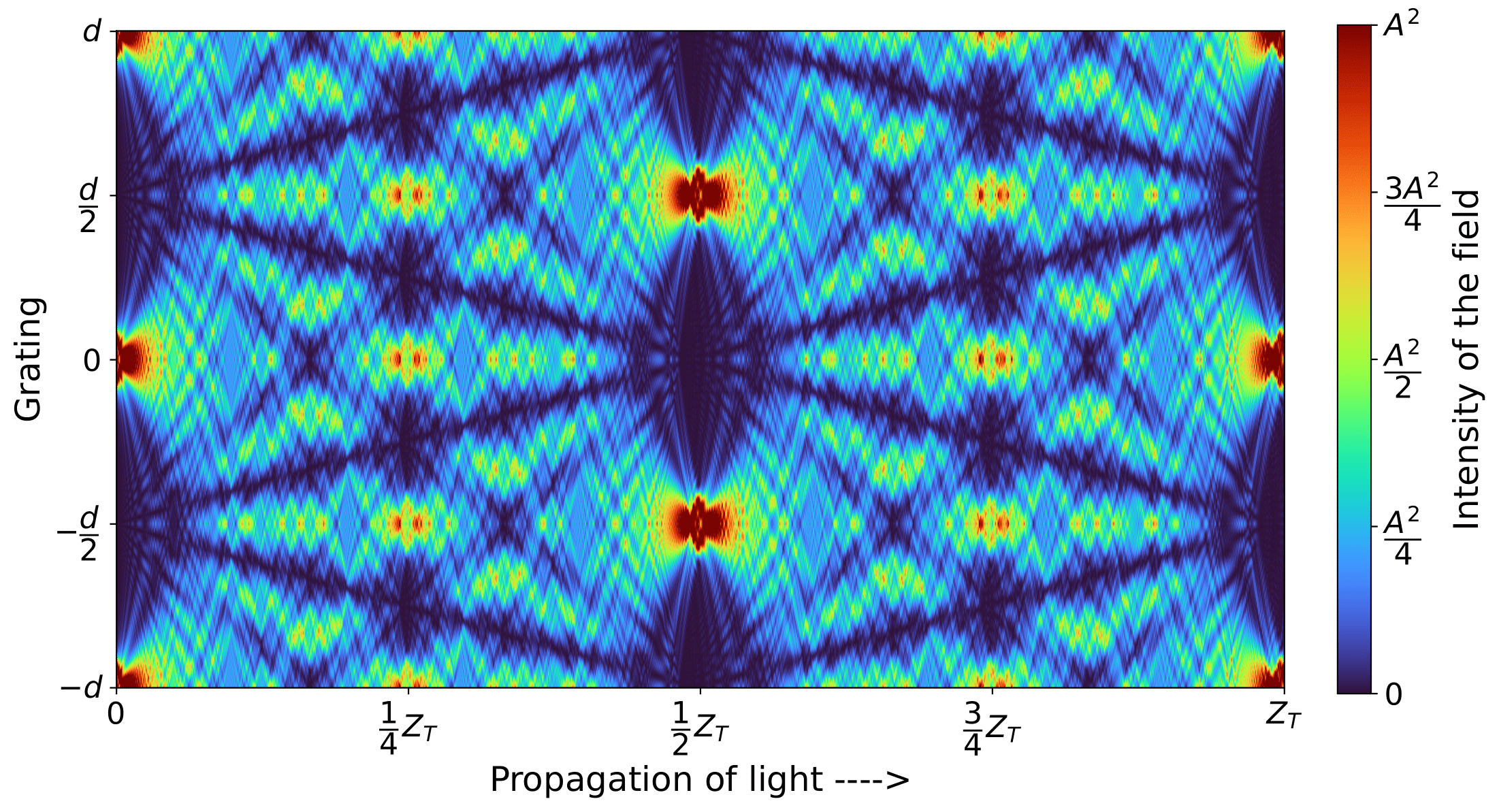}
      \caption{Carpet with $\frac{d}{\lambda}=10$ and $\frac{l}{\lambda}=2$.}
      \label{fig:d_lambda=10_w_lambda=2_energy_carpet}
    \end{subfigure}
    
    \begin{subfigure}{.45\textwidth}
      \centering
      \includegraphics[width=\linewidth]{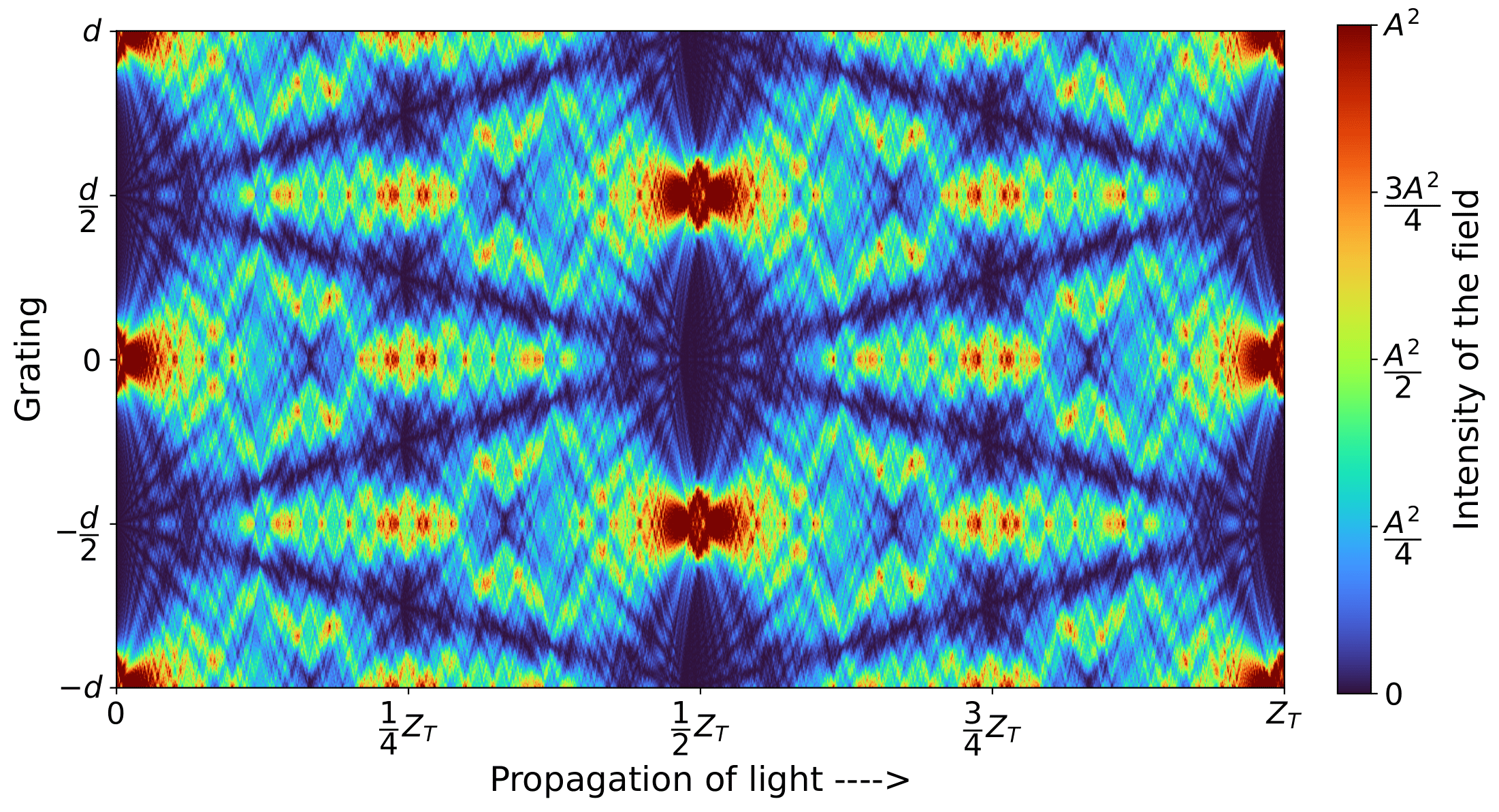}
      \caption{Carpet with $\frac{d}{\lambda}=20$ and $\frac{l}{\lambda}=5$.}
      \label{fig:d_lambda=20_w_lambda=5_energy_carpet}
    \end{subfigure}
    \begin{subfigure}{.45\textwidth}
      \centering
      \includegraphics[width=\linewidth]{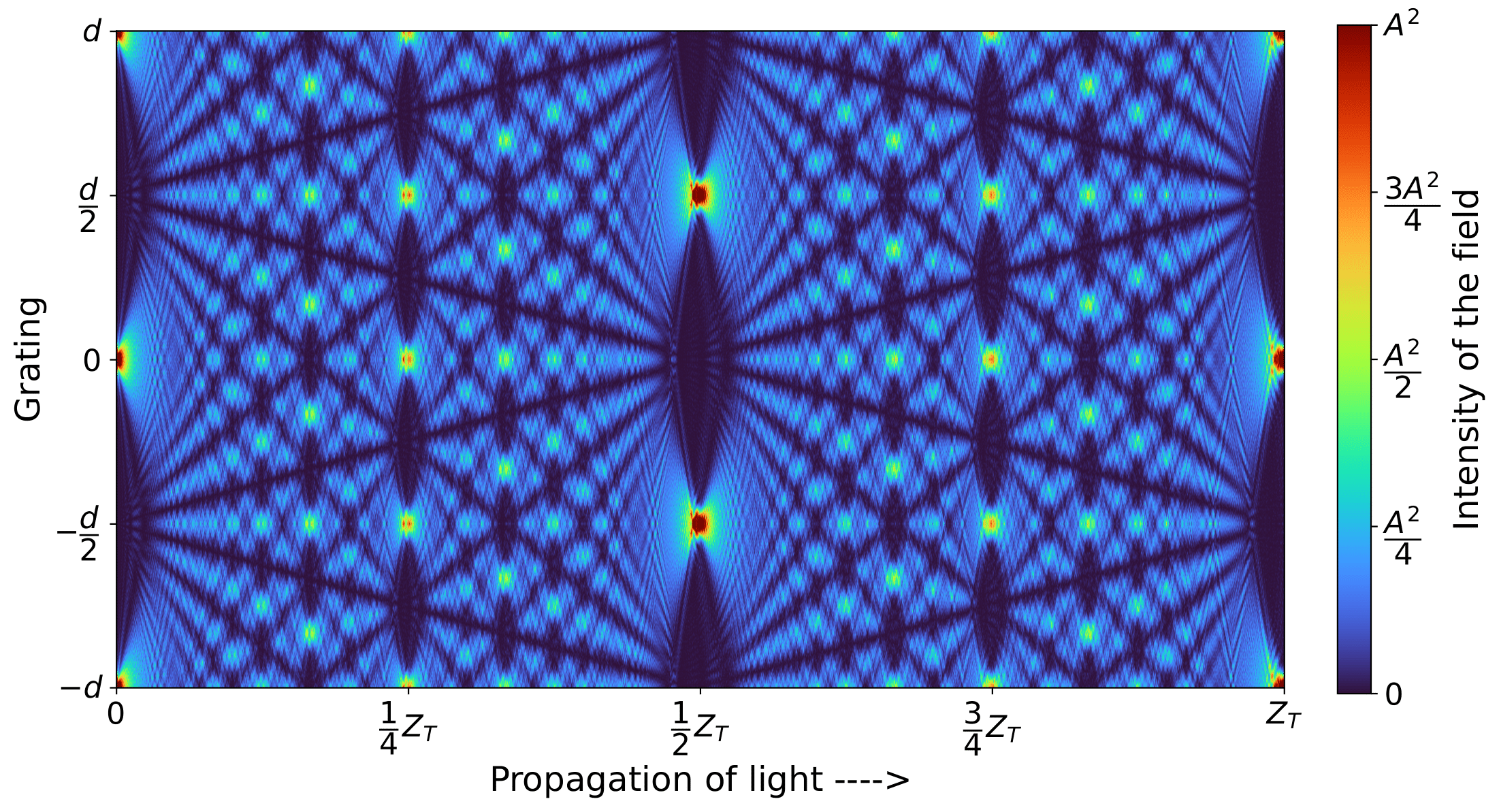}
      \caption{Carpet with $\frac{d}{\lambda}=20$ and $\frac{l}{\lambda}=2$.}
      \label{fig:d_lambda=20_w_lambda=2_energy_carpet}
    \end{subfigure}

    \begin{subfigure}{.45\textwidth}
      \centering
      \includegraphics[width=\linewidth]{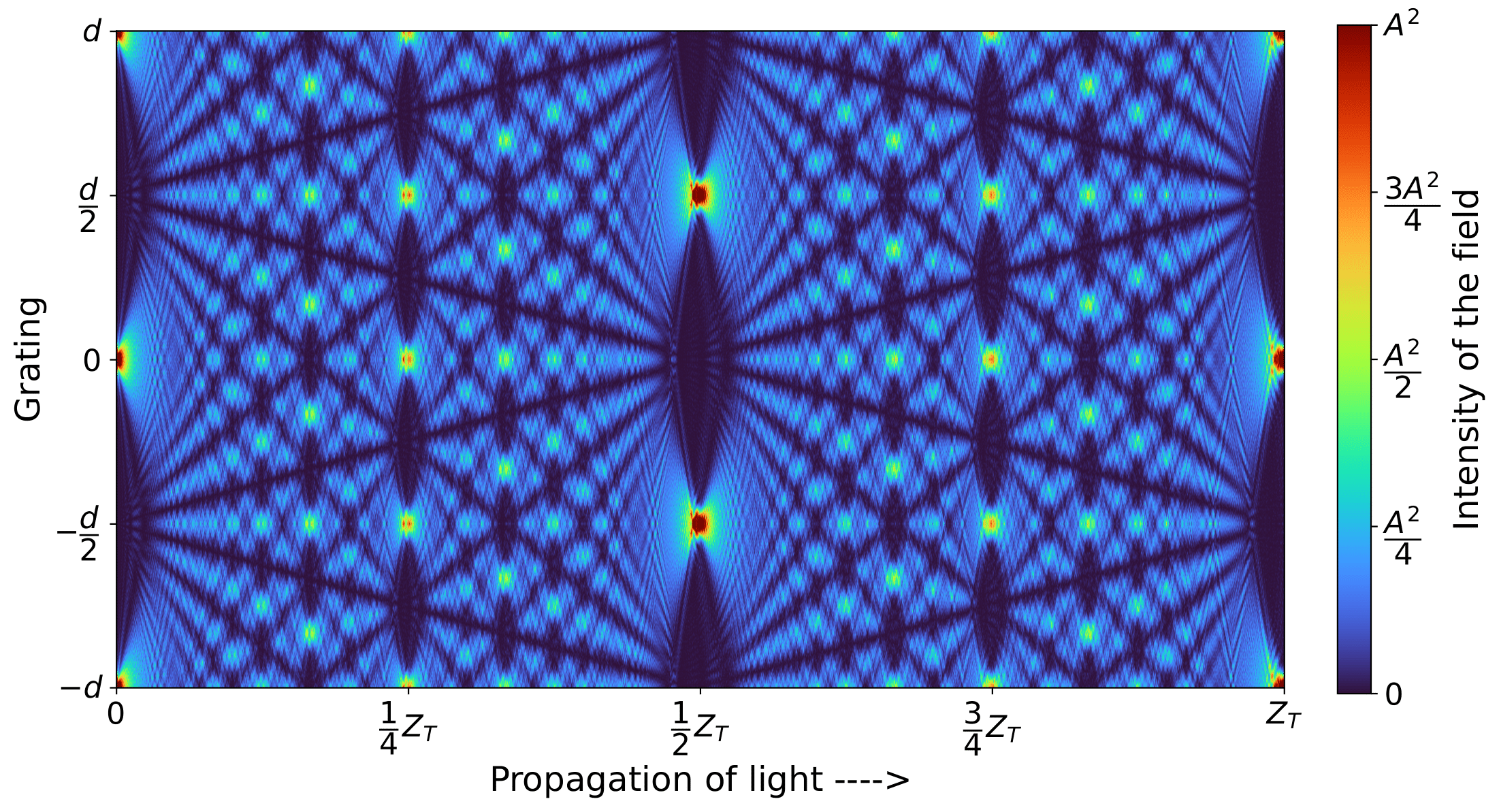}
      \caption{Carpet with $\frac{d}{\lambda}=50$ and $\frac{l}{\lambda}=5$.}
      \label{fig:d_lambda=50_w_lambda=5_energy_carpet}
    \end{subfigure}
    \begin{subfigure}{.45\textwidth}
      \centering
      \includegraphics[width=\linewidth]{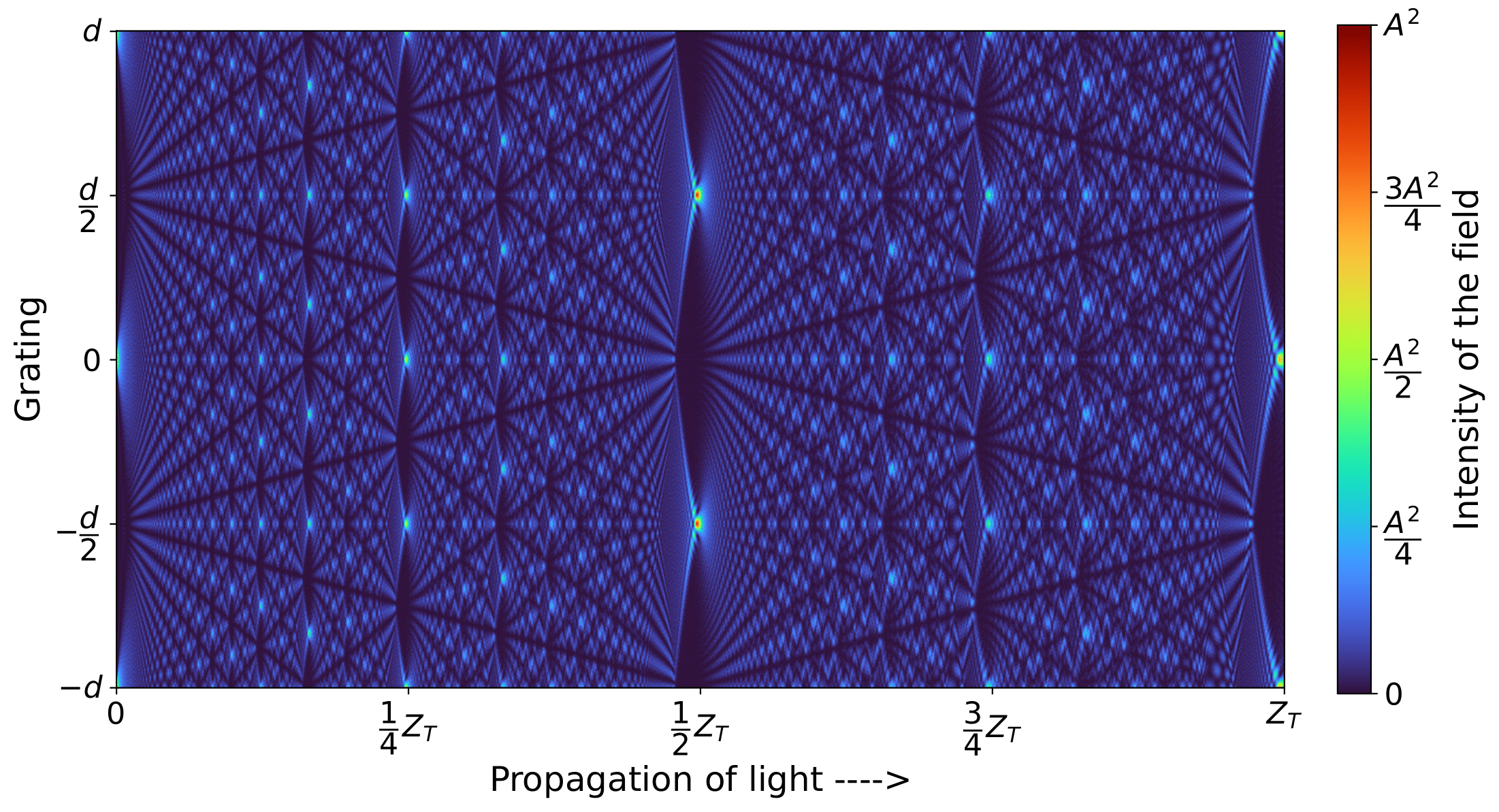}
      \caption{Carpet with $\frac{d}{\lambda}=50$ and $\frac{l}{\lambda}=2$.}
      \label{fig:d_lambda=50_w_lambda=2_energy_carpet}
    \end{subfigure}

    \begin{subfigure}{.45\textwidth}
      \centering
      \includegraphics[width=\linewidth]{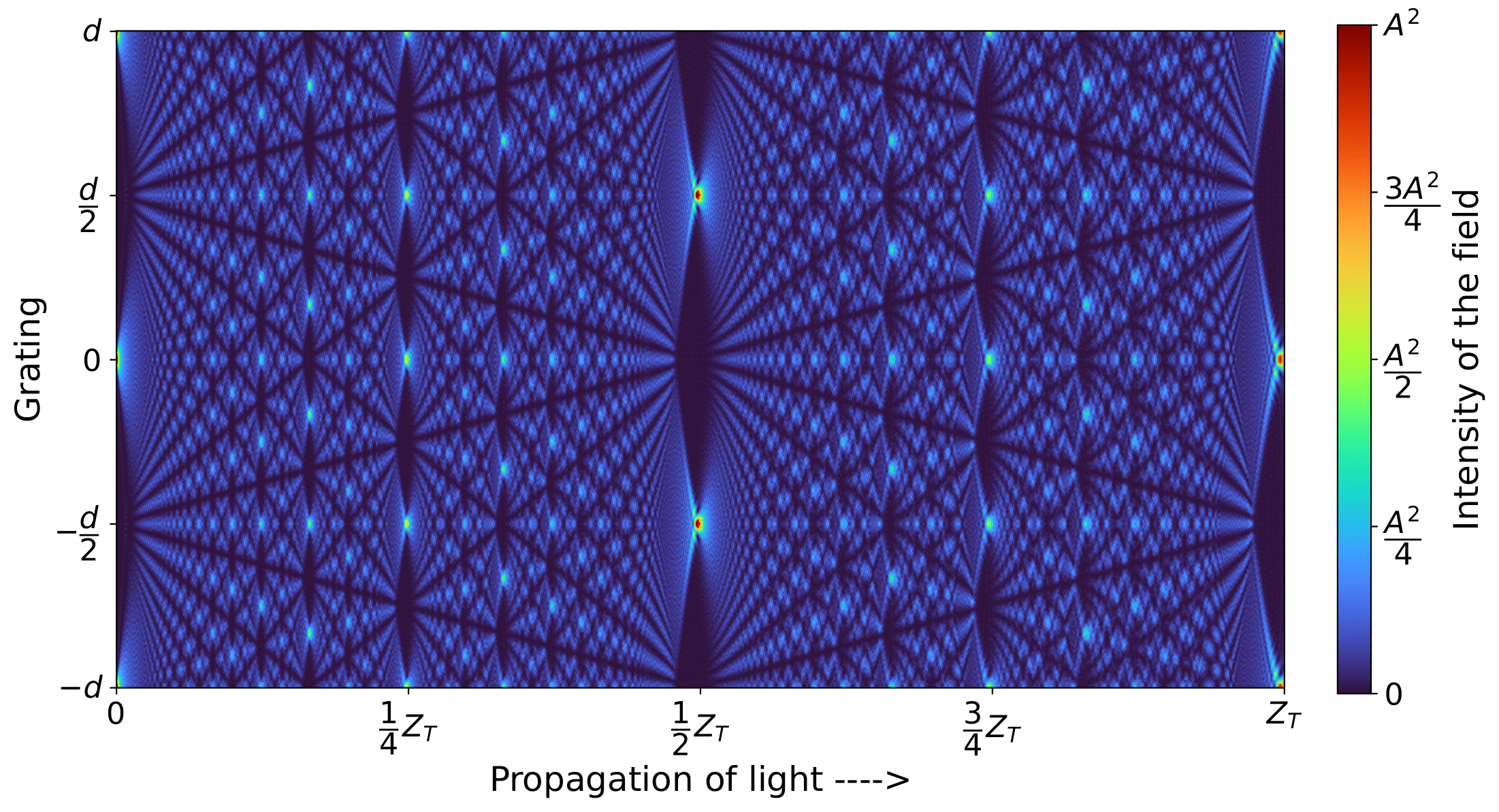}
      \caption{Carpet with $\frac{d}{\lambda}=100$ and $\frac{l}{\lambda}=5$.}
      \label{fig:d_lambda=100_w_lambda=5_energy_carpet}
    \end{subfigure}
    \begin{subfigure}{.45\textwidth}
      \centering
      \includegraphics[width=\linewidth]{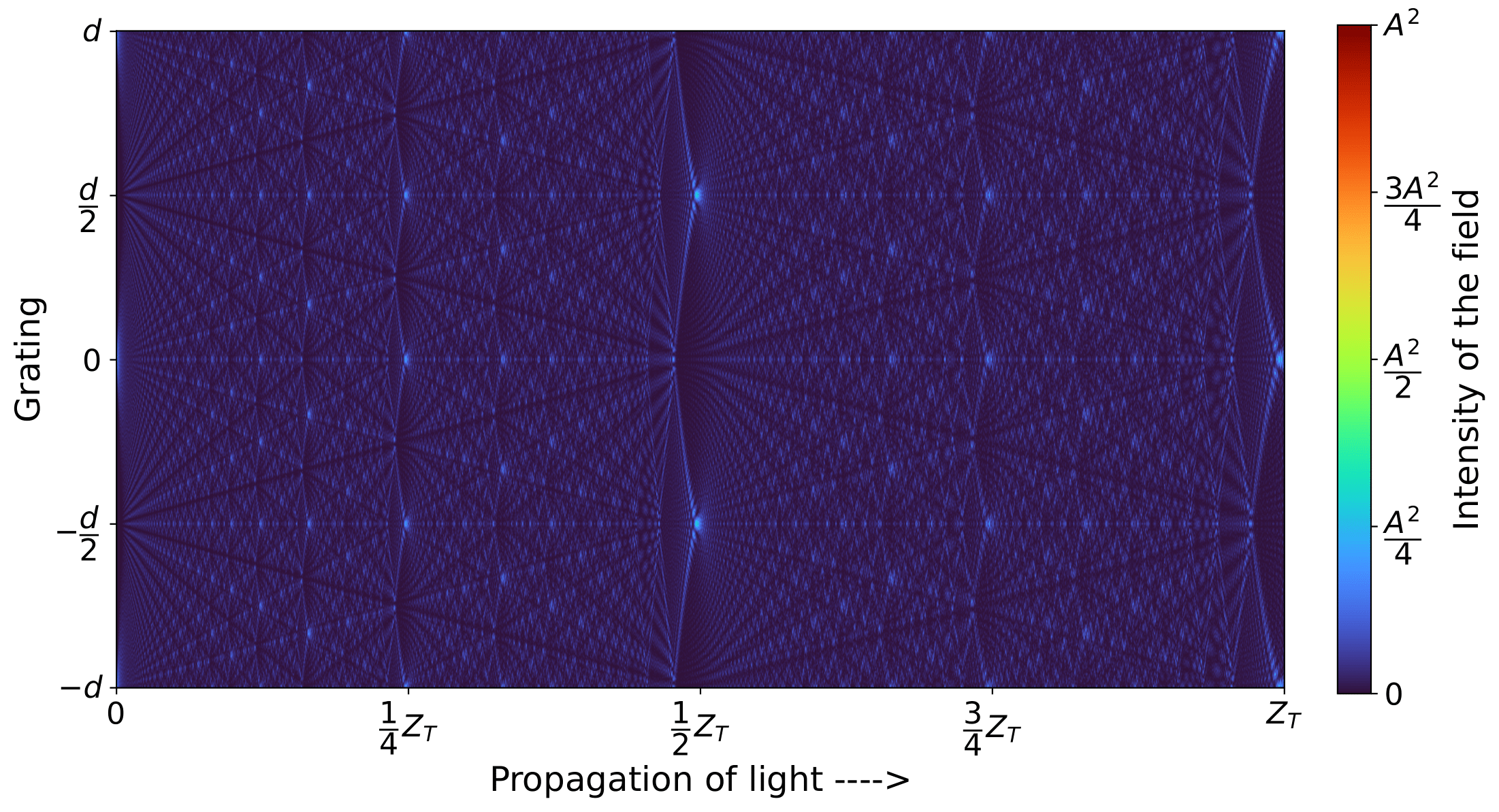}
      \caption{Carpet with $\frac{d}{\lambda}=100$ and $\frac{l}{\lambda}=2$.}
      \label{fig:d_lambda=100_w_lambda=2_energy_carpet}
    \end{subfigure}
    
    \caption{Multiple Talbot carpets obtained by squaring and truncating the envelope of the intensity $\abs{U}$, from equation \eqref{eq:U_definition} at the $N = 5\,\frac{d}{\lambda}$-th term, as detailed in Section \ref{section:graphing_carpets}. A warmer colour means higher intensity. Figure \ref{fig:d_lambda=5_w_lambda=5_energy_carpet} is just a propagating wave in free space.}
    \label{fig:multiple_energy_carpets_transient}\vspace{-4mm}
\end{figure}

We can now plot the intensity carpet $\abs{U}^2$, which we do in Figure \ref{fig:multiple_energy_carpets_transient} to see what we can learn from this new function. We notice again that the Talbot subimages multiply and become more defined as $\lambda/d, w/d \to 0$. Furthermore, in Figure \ref{fig:multiple_carpets_transient}, it can be seen that at some rational multiples of half the Talbot distance verifying $z/\frac{z_T}{2} = p/q$ with $p$ and $q$ coprime, $q$ equally spaced Talbot images appear. Again, this statement gains accuracy as the choice of parameters $l/d$ and $\lambda/d$ approaches zero. The best example of this is Figure \ref{fig:d_lambda=100_w_lambda=2_energy_carpet}.\smallskip

Additionally, we can see that at $z/\frac{z_T}{2} = p/q$ each subimage's amplitude is close to $A/\sqrt{q}$ --which is necessary for energy conservation over the z-direction. This also gains accuracy as the contribution of the evanescent waves decreases, as these are waves that \emph{dissipate} energy along the $z$ direction.\smallskip

Finally, we notice that there are some lines of darkness. Indeed, we can see that the straight path between $x=d/2$, $z=0$ and $x=0$, $z=z_T/2$ is dark. More generally, any path from $x=d/2$, $z=0$ to $x=\nu d$, $z=z_T/2$ with $\nu\in\mathbb{Z}$ is dark. Later, in Subsection \ref{section:dark_paths}, we will use this as a test for whether our limit is correct.

\section{The paraxial limit}\label{section:paraxial}

\subsection{The rational Talbot effect}

We have seen that as the ratio $\lambda/d$ diminishes, the number of subimages increases. Motivated by this, we will investigate what happens when $\lambda/d \to 0$. This is called the \emph{paraxial limit}, for which we expect to have a subimage at every rational distance $z/(z_T/2) = p/q$, where $p,q$ are coprime integers. Let us show that this is the case.\smallskip

We first define the \emph{reduced variables} $\xi$ and $\zeta$ as
\begin{align}
    \xi &:= \dfrac{x}{d} & \zeta &:= \dfrac{z}{z_T/2} = \dfrac{\lambda}{d^2} z
\end{align}

In the paraxial limit we have that $\omega \gg k_n$, so,
\begin{equation*}
    -z \sqrt{\omega^2 - k_n^2} \approx - z \omega + z\dfrac{k_n^2}{\omega} = - z \omega + \pi \zeta n^2,
\end{equation*}
where we have used the fact that $z_T = 2d^2/\lambda = d^2\omega/\pi $. Furthermore, in this limit the real exponential terms do not appear, so equation \eqref{eq:U_definition} reduces to
\begin{equation*}
    U_{\text{Par}} (\xi,\zeta) = e^{-i\omega z} \sum_n \hat{g}_n e^{i 2\pi \xi n + i \pi \zeta n^2}.
\end{equation*}
We also remark that this object solves the free Schrödinger equation 
\begin{equation*}
    -i\,\partial_\zeta U_{\text{Par}} = - \left(\partial_{\xi\xi} + \dfrac{\omega z_T}{2}\right) U_{\text{Par}},
\end{equation*}
which shows that a Talbot effect is to be expected in the context of quantum mechanics. In this case the propagation direction $z$ is playing the role of time, so we have \emph{spacial quantum revivals}, instead of \emph{temporal} ones. A by-product of this is that we find that the positions and momenta are discrete for rational times! Further on this topic can be found in \cite{Eceizabarrena2021}.\smallskip

Let us show that there are $q$ perfect Talbot subimages in the $\zeta = \nu + p/q$ line, with $p,q$ two coprime integers and $\nu$ another integer. For this we do the decomposition
\begin{equation*}
    e^{i\pi p n} e^{i\pi p/q\, n^2} = \sum_{r=0}^{q-1} c_r e^{2\pi i rn/q}.
\end{equation*}
In order to determine the coefficients $c_m$ we use the completeness relation
\begin{equation*}
    \sum_{n=0}^{q-1} e^{2\pi i a n/q } e^{-2\pi i b n/q} = q \delta_{a,b},
\end{equation*}
so that
\begin{equation*}
    \sum_{n=0}^{q-1} e^{i\pi \frac{n^2 p + pq n}{q} } e^{-2\pi i a n/q} = \sum_{n=0}^{q-1} \sum_{r=0}^{q-1} c_r e^{2\pi i rn/q} e^{-2\pi i a n/q} = q c_a,
\end{equation*}
thus the coefficients are
\begin{equation*}
    c_r = \dfrac{1}{q}\sum_{n=0}^{q-1} e^{2\pi i \frac{p/2 n^2 + (pq/2 - r) n}{q}} \equiv \dfrac{1}{q} \mathbb{G}\left(p/2, pq/2 - r , q \right),
\end{equation*}
where we have defined the generalised Gauss sum $\mathbb{G}(p, r, q) := \sum_{n=0}^{q-1} e^{2\pi i \frac{n^2 p + r n}{q}}$. Some properties of these sums are explained in Appendix \ref{appendix:GQGS}.\smallskip

We now insert this into $U$ to get
\begin{equation*}
    U_{\text{Par}} (\xi, \nu + p/q) = e^{-i\omega z}  \sum_{r=0}^{q-1} c_r \sum_n \hat{g}_n e^{i 2\pi \xi n + i\pi\nu n + i 2 \pi r n/q} = \sum_{r=0}^{q-1} c_r U_{\text{Par}}\left(\xi + \dfrac{r}{q} + \dfrac{\nu + p}{2},\zeta=0 \right).
\end{equation*}
We thus see that we get $q$ exact replicas of the distribution at $z=0$ of amplitude $\abs{c_r}$, but shifted in $\xi$ by $(\nu + p)/2 + r/q$. Of course, in a real scenario the replicas will not be exact due to the non-paraxial terms. Furthermore, we can use the fact that $\abs{\mathbb{G}\left(p/2, pq/2 - r , q \right)} = 1/\sqrt{q}$, which is proven in Appendix \ref{appendix:GQGS} to get that
\begin{equation}
    \abs{U_{\text{Par}} (\xi, \nu + p/q)}  = \dfrac{1}{\sqrt{q}} \sum_{r=0}^{q-1} U_{\text{Par}} \left(\xi + \dfrac{r}{q} + \dfrac{\nu + p}{2},\zeta=0 \right),
\end{equation}
from which we extract that each of the subimages has the same amplitude.

\subsection{Convergence of the solution}

In the last section we have established that it is plausible that in the paraxial limit $\lambda/d \to 0$ there are $q$ perfect subimages at distances $\zeta=p/q$, but we have not proven it. Due to the oscillatory nature of the sum, we will only prove a convergence in the $L^2$ sense. Physically, this means that want the \emph{mean energy over one period in the $x$-direction} of the non-paraxial terms to vanish. Indeed, what we want to show is that

\begin{theo}
    Given the parameter $\epsilon := \lambda/d$, the solution to the Helmholtz equation
    \begin{equation*}
        U_\epsilon (\xi, \zeta) = \sum_n \hat{g}_n \times e^{2\pi i \xi n} \begin{cases}
           \exp \left(- i 2\pi \zeta/\epsilon^2 \sqrt{1 - (n\epsilon)^2} \right), &  0< n < 1/\epsilon \\ 
           \exp \left(- 2\pi \zeta/\epsilon^2\sqrt{(n\epsilon)^2 - 1}\right), &  1/\epsilon < n <\infty
        \end{cases},
    \end{equation*}
    and the ideal paraxial solution
    \begin{equation*}
        U_{\text{Par}} (\xi,\zeta) = e^{-i\omega z} \sum_n \hat{g}_n e^{i 2\pi \xi n + i \pi \zeta n^2},
    \end{equation*}
    then, for any grating $g(x)\in L^2$ with Fourier coefficients $\hat{g}_n$, $U$ converges to $U_{\text{Par}}$ in the $L^2$ sense
    \begin{equation}
        \lim_{\epsilon \to 0} \int_0^1 \abs{U_\epsilon (\xi, \zeta) - U_{\text{Par}} (\xi, \zeta)}^2 \,d\,\xi = 0.
    \end{equation}
\end{theo}\smallskip

\emph{Proof:} For this we split $U$ into two distinct sums:
\begin{align*}
    U^{\text{I}}_\epsilon (\xi, \zeta)	&:= \sum_{\abs{n} \leq \lfloor 1/\epsilon \rfloor} \hat{g}_n e^{- i 2\pi \zeta/\epsilon^2 \sqrt{1 - (n\epsilon)^2}}\, e^{ 2\pi i \xi n} ,\\
    U^{\text{II}}_\epsilon (\xi, \zeta)	&:= \sum_{\abs{n} > \lfloor 1/\epsilon \rfloor} \hat{g}_n e^{- 2\pi \zeta/\epsilon^2 \sqrt{(n\epsilon)^2 - 1}}\, e^{ 2\pi i \xi n}.
\end{align*}
Our task is then to bound
\begin{equation*}
    \lim_{\epsilon \to 0} \int_0^1 \abs{U_\epsilon - U_{\text{Par}}}^2  \,d\,\xi \leq \lim_{\epsilon \to 0} \int_0^1 \abs{U^{\text{I}}_\epsilon - U_{\text{Par}}}^2  \,d\,\xi + \lim_{\epsilon \to 0} \int_0^1 \abs{U^{\text{II}}_\epsilon}^2  \,d\,\xi,
\end{equation*}
where we have used the orthogonality of the complex exponentials. We mow show that both these integrals are zero.\smallskip

\emph{ Bounding $\abs{U^{\text{I}}_\epsilon - U_{\text{Par}}}$:} First of all we expand the square root using its Taylor series
\begin{equation*}
    z \sqrt{\omega^2 - k_n^2} = \omega z\,\sqrt{1 - (n\epsilon)^2} \approx \omega z - \pi \zeta n^2 + O \left( \epsilon^2 \, n^4 \right),
\end{equation*}
thus we have that
\begin{equation*}
U^{\text{I}}_\epsilon (\xi, \zeta) - U_{\text{Par}} (\xi, \zeta) = e^{- i \omega z} \sum_{\abs{n} \leq \lfloor 1/\epsilon \rfloor} \hat{g}_n \, e^{2\pi i n \xi + \pi i \zeta n^2} \left( e^{i O \left( \epsilon^2 \, n^4 \right)} - 1 \right),
\end{equation*}
so that
\begin{equation*}
    \int_0^1 \abs{U^{\text{I}}_\epsilon (\xi, \zeta) - U_{\text{Par}} (\xi, \zeta)}^2  \,d\,\xi = \sum_{\abs{n} \leq \lfloor 1/\epsilon \rfloor} \abs{\hat{g}_n}^2  \abs{ e^{i O \left( \epsilon^2 \, n^4 \right)} - 1 }^2,
\end{equation*}
where all the cross terms are zero because of orthogonality. Now, we also have that $ \abs{ e^{i O \left( \epsilon^2 \, n^4 \right)} - 1 }^2$ is upper bounded, $ \abs{ e^{i O \left( \epsilon^2 \, n^4 \right)} - 1 }^2 < 4$, and converges pointwise to zero as $\epsilon\to0$. Furthermore, as $g(x) \in L^2$,
\begin{equation*}
    \sum_{n=1}^\infty \abs{\hat{g}_n}^2 \abs{ e^{i O \left( \epsilon^2 \, n^4 \right)} - 1 }^2 \leq \sum_{n=1}^\infty 4 \abs{\hat{g}_n}^2,
\end{equation*}
which is convergent. Thus, by Weierstrass' M-test the series converges uniformly and we have that
\begin{equation*}
    \lim_{\epsilon \to 0}\sum_{\abs{n} \leq \lfloor 1/\epsilon \rfloor} \abs{\hat{g}_n}^2  \abs{ e^{i O \left( \epsilon^2 \, n^4 \right)} - 1 }^2 = \sum_{\abs{n}}^\infty \lim_{\epsilon \to 0} \abs{\hat{g}_n}^2  \abs{ e^{i O \left( \epsilon^2 \, n^4 \right)} - 1 }^2 = 0,
\end{equation*}
so
\begin{equation}
    \lim_{\epsilon \to 0} \int_0^1 \abs{U^{\text{I}}_\epsilon (\xi, \zeta) - U_{\text{Par}} (\xi, \zeta)}^2  \,d\,\xi = 0.
\end{equation}\smallskip

\emph{Bounding $\abs{U^{\text{II}}_\epsilon}$:} We proceed as before, and find that
\begin{equation*}
    \int_0^{1} \abs{U^{\text{II}}_\epsilon (\xi, \zeta)}^2 \, d\xi = \int_0^{1}\,\abs{ \sum_{\abs{n}>\lfloor 1/\epsilon \rfloor} \hat{g}_n \,e^{- 2\pi \zeta/\epsilon^2 \sqrt{(n\epsilon)^2 - 1}}\,e^{2\pi i n \xi} }^2 \, d\xi \lesssim\, \sum_{\abs{n}>\lfloor 1/\epsilon \rfloor} \abs{\hat{g}_n}^2.
\end{equation*}
But this is just the tail of a $L^2$-convergent series, so
\begin{equation}
    \lim_{\epsilon \to 0} \int_0^{1} \abs{U^{\text{II}}_\epsilon (\xi, \zeta)}^2 \, d\xi \leq \lim_{\epsilon \to 0} \sum_{\abs{n}>\lfloor 1/\epsilon \rfloor} \abs{\hat{g}_n}^2 = 0 .
\end{equation}
\hfill $\square$

\section{Dirac delta source and the fractal Talbot effect}\label{section:ideal}
\subsection{Ideal form of the solution}

Building upon last section's ideas and what we have seen in Figure \ref{fig:multiple_energy_carpets_transient}, we want to study what shape the carpet takes in the paraxial limit when we also have a very thin grating. That is, the limit with $\epsilon = \lambda/d \to 0$ and $\delta := l/d\to 0$. In this limit the Fourier coefficients of the Ronchi grating are $\hat{g}_n = \frac{d}{l} \frac{\sin{\frac{n\pi l}{d}}}{n \pi} \approx \frac{d}{l} \frac{\frac{n\pi l}{d}}{n \pi} = 1$. Thus, in such limit, the Ronchi grating becomes a Dirac comb grating. Thus
\begin{align}
\begin{split}
    U_{\text{ideal}}\left(\xi, \zeta = \nu + \dfrac{p}{q}\right) &= e^{-i\omega z}\,\sum_{r = 0}^{q-1} c_r\, U\left( \xi + \dfrac{r}{q} + \dfrac{\nu + p}{2} \right),\\
    &= e^{-i\omega z}\, \sum_{m\in\mathbb{Z}} \sum_{r = 0}^{q-1} c_r\, \delta\left( \xi + \dfrac{r}{q} + \dfrac{\nu + p}{2} - m\right),\\
    &= e^{-i\omega z}\, \sum_{m\in\mathbb{Z}} c_m\, \delta\left( \xi - \dfrac{m}{q} + \dfrac{\nu + p}{2} \right),\label{eq:U_ideal_deltas}
\end{split}
\end{align}
where $c_m = \frac{1}{q}\sum_{n=0}^{q-1} e^{2\pi i \frac{p/2 \, n^2 + (pq/2 - m) n}{q}} \equiv \dfrac{1}{q} \mathbb{G}\left(p/2, pq/2 - m , q \right)$. Now, we have that $\abs{\mathbb{G}\left(p/2, m , q \right)} = \sqrt{q}$. Thus
\begin{equation*}
    \abs{U_{\text{ideal}} \left(\xi, \zeta = \nu + \dfrac{p}{q}\right)} = \dfrac{1}{\sqrt{q}} \sum_{m\in\mathbb{Z}}  \delta\left( \xi - \dfrac{m}{q} + \dfrac{\nu + p}{2} \right).
\end{equation*}
From here we can see that each of the sub images has the same intensity, and that energy is conserved along $x$. Furthermore, we can see in Figure \ref{fig:d_lambda=100_w_lambda=5_energy_carpet} that a carpet very similar to what we would get from this last equation.\smallskip

Nevertheless, equation \eqref{eq:U_ideal_deltas} is only valid for $\zeta = \nu + p/q$. When $\zeta$ is not a rational number, we have
\begin{equation}\label{eq:U_ideal}
    U_{\text{ideal}} \left(\xi, \zeta \right) = e^{-i\omega z} \sum_n e^{i 2\pi \xi n + i \pi \zeta n^2},
\end{equation}
and it was shown that this object has a fractal structure for $\zeta$ irrational \cite{Berry1996, Kapitanski1999, Rodnianski2000}.\smallskip

In this case, the solution to the Helmholtz equation $U_\epsilon$ does not converge to $U_{\text{ideal}}$ in $L^2$. Indeed the many phase cancellations make the sum very delicate, so we can only prove a weaker convergence. Eceizabarrena showed in \cite{Eceizabarrena2021} that a convergence can be proven for $\hat{g}_n = 1$, an exact Dirac comb, but not in $L^2$. Instead, he proved that
\begin{equation*}
    \lim_{\epsilon\to 0} U_\epsilon = U_{\text{ideal}},\quad \text{in } H^{-s} \text{ with } s > \dfrac{1}{2}, 
\end{equation*}
where $H^{-s}$ is a Sobolev space.

\subsection{Dark paths}\label{section:dark_paths}

Another property that we can extract from this ideal distribution is the existence of the infinitely many dark lines in the Talbot carpet, as mentioned above. Indeed, we have seen that these are lines that go from $x=d/2$, $z=0$ to $x=(\nu + 1/2) d$, $z=z_T$ with $\nu\in\mathbb{Z}$, so we may parametrise them as
\begin{align}\label{eq:parameterisation_dark_lines}
    \xi (t) &= \dfrac{1}{2} + \nu t,      &        \quad \zeta (t) &= 2t.
\end{align}
Thus, at a rational time $\zeta = 2t = 2p/q$, so 
\begin{equation*}
    \xi - \dfrac{m}{q} + \dfrac{2p}{2} = \dfrac{1}{2} + \nu \dfrac{p}{q} - \dfrac{m}{q} + p = \dfrac{q + 2\nu p - 2m + 2pq}{2q}.
\end{equation*}
We can see that this can only be zero if $q$ is even. But this is impossible if we require $2p$ and $q$ to be coprime, which is required to use \eqref{eq:U_ideal_deltas}. This means that none of the arguments of the Dirac deltas in $U_{\text{Ideal}}$ will become zero, thus $U_{\text{Ideal}} = 0$ for the lines parametrised by \eqref{eq:parameterisation_dark_lines}. We thus have that there will be no light in such lines at rational times.

\section{Summary}

In this paper we have shown that it is possible to derive the optical Talbot effect from classical electrodynamics. Indeed, by considering the wave equation for the electromagnetic four-potential under the appropriate boundary conditions we have established a well-posed PDE problem in \eqref{pde_statement} and showed that its solution is \eqref{eq:transient_solution}. The advantage of this approach, apart from being better motivated than scalar theory, is that it has ridden the theory of the inconsistencies of having to deal with the Helmholtz equation and the Kirchhoff boundary conditions. On top of this, we derived an equation with a temporal dependency. Finally, we are able to give a precise physical interpretation to the scalar field: $A^y$.\smallskip

We proved that our solution $u$ is consistent with the usual expressions that we can find in the literature \eqref{eq:stationary_solution}, as $u$ asymptotically evolves into $u_{\text{Helm}}$. We have also written a code to plot the time evolution of the field and the formation of the Talbot effect, which can be found in \cite{gabriel_ybarra_marcaida_pytalbot_2025}. This code might be use to explore the Talbot pattern created by other gratings.\smallskip

We used the simulations to build some intuition for what shape the field should have in the paraxial limit, and we proposed an ideal distribution for such a field. We show that such a distribution would exhibit the rational Talbot effect. That is, at every plane with $\zeta = p/q$ we will have $q$ equally spaced subimages. We prove that $U$ converges to this ideal paraxial distribution in the $L^2$ sense.

\section*{Acknowledgements}

I am deeply grateful to Luis Vega for presenting this topic to me several years ago, for his many insightful comments and for his invaluable support. I also want to thank Miguel Fernández de Retana for his much appreciated help with the pyTalbot code and Paula Heim for her detailed suggestions. This work was developed during my stays at the Basque Centre for Applied Mathematics.

%%%%%%%%%%%%%%%%%%%%%%% BIBLIOGRAFÍA
\printbibliography[heading=bibintoc, title={References}]
\newpage

%%%%%%%%%%%%%%%%%%%%%%% APENDIXES
\appendix
\addtocontents{toc}{\protect\setcounter{tocdepth}{-1}}
\section{\texorpdfstring{Computing the inverse Laplace transform of $e^{-z\sqrt{s^2+k^2}}$}{Computing the inverse Laplace transform of Exp(-z Sqrt(s2+k2))}}\label{appendix:laplace}

Let us prove \eqref{eq:invLT_G}, which states that
\begin{equation*}
    \mathcal{L}_t^{-1}\{ F(s)\}(s) = \mathcal{L}_t^{-1}\left\{ e^{-z\sqrt{k^2+s^2}} \right\}(s) = \delta(t-z) - \theta(t-z) k_n z \dfrac{J_1\left( k_n\sqrt{t^2-z^2} \right)}{\sqrt{t^2-z^2}}.
\end{equation*}
We start by defining the functions
\begin{equation*}
    G(s) := F(s) - e^{-zs},
\end{equation*}
and
\begin{equation*}
    H(s) := \dfrac{d}{ds} [e^{zs}\,G(s)] = \dfrac{d}{ds} \left[e^{zs-z\sqrt{k^2+s^2}} + 1\right] = z \dfrac{\sqrt{k^2+s^2} - s}{\sqrt{k^2+s^2}} e^{zs-z\sqrt{k^2+s^2}}.
\end{equation*}
We perform a Taylor expansion of this
\begin{equation*}
    H(s) = z \dfrac{\sqrt{k^2+s^2} - s}{\sqrt{k^2+s^2}} \sum_{n=0} \dfrac{\left( zs-z\sqrt{k^2+s^2} \right)^n}{n!} = kz \sum_{n=0} \dfrac{(-kz)^n}{n!} \dfrac{\left(\sqrt{k^2+s^2} - s \right)^{n+1}}{k^{n+1} \sqrt{k^2+s^2}},
\end{equation*}
and take the inverse Laplace transform of this expression, using the identity
\begin{equation*}
    \mathcal{L}_t\left\{ J_n(kt)\theta(t) \right\}(s) = \dfrac{\left(\sqrt{k^2+s^2} - s \right)^{n}}{k^{n} \sqrt{k^2+s^2}}, \quad n>-1, \text{ Re}(s)>0.
\end{equation*}
So we find that
\begin{equation*}
    h(t) \equiv \mathcal{L}_s^{-1}\left\{ H(s) \right\}(t) = kz \sum_{n=0} \dfrac{(-kz)^n}{n!} J_{n+1}(kt)\theta(t).
\end{equation*}
Now, we want to use the multiplication theorem for Bessel functions
\begin{equation*}
    \lambda ^{-\nu }J_{\nu }(\lambda x)=\sum _{n=0}^{\infty }{\frac {1}{n!}}\left({\frac {\left(1-\lambda ^{2}\right)x}{2}}\right)^{n}J_{\nu +n}(x),
\end{equation*}
to simplify the expression. We thus rewrite
\begin{align*}
\begin{split}
    \sum_{n=0} \dfrac{(-kz)^n}{n!} J_{n+1}(kt)\theta(t) &= \sum_{n=0} \dfrac{(1 - \lambda^2)^n (kt)^n}{n!} J_{n+1}(kt)\\
    &= \lambda^{-1} J_1 (\lambda kt),
\end{split}
\end{align*}
where $\lambda = \sqrt{2z/t + 1}$. So we find
\begin{equation*}
    h(t) = \dfrac{kz}{\sqrt{2z/t + 1}} J_1 \left(\sqrt{2z/t + 1} kt\right) \theta(t) = \dfrac{kz\sqrt{t}}{\sqrt{2z + t}} J_1 \left(k\sqrt{t(2z +t)}\right) \theta(t).
\end{equation*}
Now,
\begin{equation*}
    \mathcal{L}_t \{t \, g(t)\}(s) = -\dfrac{d}{ds} G(s) = -\dfrac{d}{ds} [e^{-zs} e^{zs} G(s)] = z e^{-zs} [e^{zs} G(s)] - e^{-zs} \dfrac{d}{ds}[e^{zs} G(s)] = z G(s) - e^{-zs} H(s),
\end{equation*}
so we find that
\begin{equation*}
    \mathcal{L} \{(z-t)\,g(t)\} (s) = e^{-zs} H(s),
\end{equation*}
which implies
\begin{equation*}
    g(t)\theta(t) = \dfrac{1}{z-t} \mathcal{L}_s^{-1}\left\{ e^{-zs} H(s) \right\}(t) = - \dfrac{h(t-z)}{t-z} \theta(t-z).
\end{equation*}
This gives us
\begin{equation*}
    g(t)\,\theta(t) = - \dfrac{kz}{\sqrt{t^2-z^2}} J_1 \left(k\sqrt{t^2-z^2}\right) \theta(t-z).
\end{equation*}
Finally we can recover $v$ using the fact that
\begin{equation*}
    G(s) = F(s) - e^{-zs} \implies f(t) = g(t) + \delta(t-z),
\end{equation*}
thus
\begin{empheq}[box=\fbox]{equation}\label{eq:appendix_result}
    f(t) = \delta(t-z) - \dfrac{kz}{\sqrt{t^2-z^2}} J_1 \left(k\sqrt{t^2-z^2}\right) \theta(t-z)
\end{empheq}
as claimed. This result can also be found in \cite{bateman1954tables}.

%%%%%%%%%%%%% APENDICE ERROR
\newpage
\section{Bounding the error in the asymptotic regime}\label{appendix:bounding_E}

We now want to see how fast the solution tends to the stationary case. For this we will study the asymptotic behaviour of the error, 
\begin{equation*}
    E_n(t,z) := k_n z \int_t^{\infty} \sin\omega(t-\tau) \dfrac{J_1\left( k_n\sqrt{\tau^2-z^2} \right)}{\sqrt{\tau^2-z^2}} \,d\tau,
\end{equation*}
where $k_n = \frac{2\pi}{d} n$, in the regime $z\ll t$, which we rewrite as
\begin{equation*}
    E_n(t,z) := k_n z \int_{k_n t}^{\infty} \sin\omega(t-u/k_n) \dfrac{J_1\left(\sqrt{u^2 - (k_n z)^2} \right)}{\sqrt{u^2 - (k_n z)^2}} \,du.
\end{equation*}
We will show using a case by case scenario that at worst
\begin{empheq}[box=\fbox]{equation}
     E_{n} = \mathcal{O}\left( \sqrt{\dfrac{z}{t}} \right)
\end{empheq}

\subsection{The case \texorpdfstring{$k_n \not\approx \omega$}{k =/= w}}

We now use the asymptotic behaviour of the Bessel function,
\begin{equation*}
    J_{\alpha }(x) = {\sqrt {\frac {2}{\pi x}}}\left(\cos \left(x-{\frac {\alpha \pi }{2}}-{\frac {\pi }{4}}\right) - \dfrac{\sin \left(x-{\frac {\alpha \pi }{2}}-{\frac {\pi }{4}}\right)}{x} + {\mathcal {O}}\left(x^{-2}\right)\right), \quad\text{for } x \gg \alpha^2-1/4,
\end{equation*}
so,
\begin{align*}
\begin{split}
    E_n(t,z) &=- k_n z \sqrt{\frac {2}{\pi}} \int_{k_n t}^{\infty} \sin\omega(u/k_n - t) \dfrac{ \cos \left( \sqrt{u^2 - (k_n z)^2} - \frac {3 \pi }{4} \right) }{\left(u^2 - (k_n z)^2 \right)^{3/4}} \left[ 1 + {\mathcal {O}}\left(\frac{1}{\sqrt{u^2 - (k_n z)^2}}\right)\right]\,du,\\
    &= - k_n z \sqrt{\frac {2}{\pi}} \int_{k_n t}^{\infty} \sin\omega(u/k_n - t) \dfrac{ \cos \left( \sqrt{u^2 - (k_n z)^2} - \frac {3 \pi }{4} \right) }{\left(u^2 - (k_n z)^2 \right)^{3/4}} \left[ 1 + {\mathcal {O}}\left(\frac{1}{k_n t}\right)\right]\,du,\\
    &= - k_n z \sqrt{\frac {2}{\pi}} I \left[ 1 + {\mathcal {O}}\left(\frac{1}{k_n t}\right)\right],
\end{split}
\end{align*}
where we have defined
\begin{equation*}
    I := \int_{k_n t}^{\infty} \sin\omega(u/k_n - t) \dfrac{ \cos \left( \sqrt{u^2 - (k_n z)^2} - \frac {3 \pi }{4} \right) }{\left(u^2 - (k_n z)^2 \right)^{3/4}}  du
\end{equation*}
We now expand
\begin{align*}
    I &\approx \int_{k_n t}^{\infty} \sin\omega(u/k_n - t) \dfrac{ \cos \left( u - \frac {3 \pi }{4} + \mathcal{O}(k_n^2 z^2/u) \right) }{\left[u + \mathcal{O}(k_n^2 z^2/u)\right]^{3/2}}  du,\\
    &\approx \int_{k_n t}^{\infty} \sin\omega(u/k_n - t) \dfrac{ \cos \left( u - \frac {3 \pi }{4} \right) }{u^{3/2}} \left[ 1 + \mathcal{O} \left(\dfrac{k_n^2 z^2}{u} \right) \right]  du,\\
    &\approx \int_{k_n t}^{\infty} \dfrac{\sin \left[ (\omega/k_n + 1) u - \omega t - \frac {3 \pi }{4} \right] + \sin \left[  (\omega/k_n - 1) u - \omega t + \frac {3 \pi }{4} \right]  }{2 u^{3/2}}  du \left[ 1 + \mathcal{O} \left(k_n z \dfrac{z}{t} \right) \right].
\end{align*}
We thus need to estimate integrals of the form $\int_t^\infty \frac{\sin (ax+b)}{x^c}\,dx$. Using integration by parts we find that 
\begin{align*}
    \int_{kt}^\infty \frac{\sin (ax+b)}{x^c}\,dx &= -\left.\frac{\cos (ax+b)}{ax^c}\right|^\infty_{kt} - c\int_{kt}^\infty \frac{\cos (ax+b)}{ax^{c+1}}\,dx\\
    &= \frac{\cos (akt+b)}{a (kt)^c} - c\int_{kt}^\infty \frac{\cos (ax+b)}{ax^{c+1}}\,dx\\
    &= \frac{\cos (akt+b)}{a (kt)^c} \left[1 + \mathcal{O}(1/kt)\right].
\end{align*}
Thus,
\begin{align*}
\begin{split}
    I &= \sqrt{\dfrac{1}{4 k_n t}} \left[ \frac{\cos ((\omega + k_n)t - \omega t - {\frac {3 \pi }{4}})}{(\omega + k_n) t} + \frac{\cos ((\omega - k_n)t -  \omega t  +{\frac {3 \pi }{4}})}{(\omega - k_n)t} \right] \left[ 1 + \mathcal{O} \left( k_n z \dfrac{z}{t} \right) \right] \left[1 + {\mathcal {O}}\left(\dfrac{1}{k_n t}\right)\right],\\
    &= \sqrt{\dfrac{1}{4 k_n t}} \left[ \frac{\cos (k_n t - {\frac {3 \pi }{4}})}{(\omega + k_n) t} + \frac{\cos ( - k_n t  +{\frac {3 \pi }{4}})}{(\omega - k_n)t} \right] \left[ 1 + \mathcal{O} \left( k_n z \dfrac{z}{t} \right) + {\mathcal {O}}\left(\dfrac{1}{k_n t}\right)\right],\\
    &= \sqrt{\dfrac{1}{k_n t}} \frac{\cos ( k_n t  - \frac {3 \pi }{4})}{t} \dfrac{\omega}{\omega^2 - k_n^2}  \left[ 1 + \mathcal{O} \left( k_n z \dfrac{z}{t} \right) + {\mathcal {O}}\left(\dfrac{1}{k_n t}\right)\right].
\end{split}
\end{align*}
We thus find that
\begin{empheq}[box=\fbox]{equation}
    E_n(t,z) = -  \dfrac{z}{t} \sqrt{\dfrac{2 k_n}{\pi t}} \cos ( k_n t  - \frac {3 \pi }{4}) \dfrac{\omega}{\omega^2 - k_n^2}  \left[ 1 + \mathcal{O} \left( k_n z \dfrac{z}{t} \right) + {\mathcal {O}}\left(\dfrac{1}{k_n t}\right)\right] = O\left(\dfrac{z}{t}\right)
\end{empheq}
This is small except whenever $k_n^2 \approx \omega^2$, as we have a resonance for $\omega = k_n$.

\subsection{The case \texorpdfstring{$k_n = \omega$}{k = w}}
In this case we have that the integral $E_n$ can be computed without approximating the Bessel function by its asymptotic series. Indeed, using the fact that
\begin{equation*}
    \dfrac{J_1\left(\sqrt{u^2 - (k_n z)^2} \right)}{\sqrt{u^2 - (k_n z)^2}} = \dfrac{J_1\left( u + \frac{k_n^2 z^2}{u} \right)}{u + \frac{k_n^2 z^2}{u} } = \dfrac{J_1\left( u \right)}{u} \left[ 1 + \mathcal{O} \left(\frac{k_n^2 z^2}{u}  \right)\right],
\end{equation*}
we have that
\begin{align*}
    E_{\omega} &= \omega z \int_{\omega t}^{\infty} \sin(\omega t-u) \dfrac{J_1\left( u \right)}{u} \left[ 1 + \mathcal{O} \left(\frac{\omega^2 z^2}{u}  \right)\right] \,du,\\ 
    &\approx \omega z  \int_{\omega t}^{\infty} \sin(\omega t-u) \dfrac{J_1\left( u \right)}{u} \,du \left[ 1 + \mathcal{O} \left(\omega z \frac{z}{t}  \right)\right],
\end{align*}
and using the fact that 
\begin{equation*}
     \int_{\omega t}^{\infty} \sin(\omega t-u) \dfrac{J_1\left( u \right)}{u}\,du = -J_0 (\omega t),
\end{equation*}
we find that
\begin{equation*}
     E_\omega \approx - \omega z  J_0 (\omega t) \left[ 1 + \mathcal{O} \left(\omega z \frac{z}{t}  \right)\right] \approx  - \omega z  \sqrt{\dfrac{2}{\pi\omega t}} \left[ 1 + \mathcal{O} \left(\omega z \frac{z}{t}  \right) + \mathcal{O} \left(\frac{1}{\omega t}  \right)\right],
\end{equation*}
Or 
\begin{empheq}[box=\fbox]{equation}
     E_\omega \approx  - \sqrt{\dfrac{z}{t}}  \sqrt{\dfrac{2 \omega z}{\pi}} \left[ 1 + \mathcal{O} \left(\omega z \frac{z}{t}  \right) + \mathcal{O} \left(\frac{1}{\omega t}  \right)\right]  = O\left(\sqrt{\dfrac{z}{t}}\right)
\end{empheq}

\subsection{The case \texorpdfstring{$k_n \approx \omega$}{k ~ w}}

We now consider the case where $k_n = \omega (1 + \delta)$, with $\delta\ll 1$. In this case we have
\begin{align*}
     E_{\omega(1 + \delta)} &= \omega z \int_{\omega (1 + \delta) t}^{\infty} \sin(\omega t - (1+\delta)u) \dfrac{J_1\left( u \right)}{u} \left[ 1 + \mathcal{O} \left((1 + \delta)^2\frac{\omega^2 z^2}{u}  \right)\right] \,du,\\
     &= \omega z \int_{\omega (1 + \delta) t}^{\infty} \left[ \sin(\omega t - u) \cos (\delta u) - \cos(\omega t - u) \sin (\delta u) \right] \dfrac{J_1\left( u \right)}{u} \left[ 1 + \mathcal{O} \left(\frac{\omega^2 z^2}{u}  \right)\right] \,du.
\end{align*}
As $\abs{\cos x} \leq 1$ and $\abs{\sin x} \leq 1$ we will clearly have that
\begin{align*}
     \abs{E_{\omega(1 + \delta)}} &\leq \omega z \left[\abs{\int_{\omega (1 + \delta) t}^{\infty} \sin(\omega t - u)  \dfrac{J_1\left( u \right)}{u}  \,du} + \abs{\int_{\omega (1 + \delta) t}^{\infty} \cos(\omega t - u)  \dfrac{J_1\left( u \right)}{u}  \,du}\right] \left[ 1 + \mathcal{O} \left(\omega z\frac{z}{t}  \right)\right],\\
     &\leq \omega z \left[\abs{\int_{\omega t}^{\infty} \sin(\omega t - u)  \dfrac{J_1\left( u \right)}{u}  \,du} + \abs{\int_{\omega t}^{\infty} \cos(\omega t - u)  \dfrac{J_1\left( u \right)}{u}  \,du}\right] \left[ 1 + \mathcal{O} \left( \delta \right) + \mathcal{O} \left(\omega z\frac{z}{t}  \right)\right].
\end{align*}
We now use the fact that
\begin{equation*}
     \int_{\omega t}^{\infty} \cos(\omega t-u) \dfrac{J_1\left( u \right)}{u}\,du = J_1 (\omega t),
\end{equation*}
which gives us that
\begin{align*}
     \abs{E_{\omega(1 + \delta)}} &\leq \omega z \left[\abs{J_0(\omega t)} + \abs{J_1(\omega t)}\right] \left[ 1 + \mathcal{O} \left( \delta \right) + \mathcal{O} \left(\omega z\frac{z}{t}  \right)\right],\\
     &\leq 2\omega z \sqrt{\dfrac{2}{\pi \omega t}} \left[ 1 + \mathcal{O} \left( \delta \right) + \mathcal{O} \left(\omega z\frac{z}{t}  \right)\right]
\end{align*}
so that
\begin{empheq}[box=\fbox]{equation}
     \abs{E_{\omega(1 + \delta)}} \leq 2\sqrt{\dfrac{z}{t}} \sqrt{\dfrac{2\omega z}{\pi}} \left[ 1 + \mathcal{O} \left( \delta \right) + \mathcal{O} \left(\omega z\frac{z}{t}  \right)\right] = \mathcal{O} \left(\sqrt{\dfrac{z}{t}}\right)
\end{empheq}
\hfill$\square$

%%%%%%%%%%%%% APENDICE SUMAS GAUSS
\newpage
\section{On generalized quadratic Gauss sums}\label{appendix:GQGS}

Let $p$, $r$ and $q$ be three integers with $r,q\,\in\,\mathbb{N}$. The generalised quadratic Gauss sum is defined as
\begin{equation}\label{definicion_generalizada}
    \mathbb{G}(p,r,q)=\sum_{n=0}^{q-1} e^{2 \pi i \frac{p n^2 + rn}{q}}.
\end{equation}
We shall prove that the modulus of this sum whenever $\gcd(p,q)=1$ is equal to
\begin{equation*}
    \abs{\mathbb{G}(p,r,q)} = \begin{cases}
        \sqrt{q}, 	& q \text{ is odd},\\
        \sqrt{2q}, 	&  2\mid q, \quad q \equiv r\pmod{4},\\
        0, 			& 2\mid q, \quad q \not\equiv r\pmod{4},
    \end{cases}
\end{equation*}\smallskip
\emph{Proof.} Introducing $\omega_q := e^{\frac{2\pi i}{q}}$, the absolute square of $\mathbb{G}$ is 
\begin{equation*}
    \abs{\mathbb{G}(p,r,q)}^2 = \mathbb{G}(p,r,q) \,\, \overline{\mathbb{G}(p,r,q)} = \sum_{m,n=0}^{q-1} \omega_q^{\,p\, (m^2-n^2) + r(m-n)} ,
\end{equation*}
which can be rewritten as
\begin{equation*}
    \abs{\mathbb{G}(p,r,q)}^2 =  \sum_{m,n=0}^{q-1} \omega_q^{\, p\, (m-n)(m+n) + r(m-n)} .
\end{equation*}
Now, let $l := n-m $, then $n=m+l$. As $l\in \left[ -m;q-m-1 \right] $ we obtain the following formula:
\begin{equation*}
\begin{split}
    \abs{\mathbb{G}}^2 &= \sum_{m=0}^{q-1}\sum_{l=-m}^{q-m-1} \omega_q^{\,-l\, \{ p \,(l+2m) + r \} } \\
    &= \sum_{m=0}^{q-1}\sum_{l=-m}^{q-m-1} \omega_q^{\,l\, \{ p \,(l+2m) + r \} }\\ &= \sum_{m=0}^{q-1} \left[ \sum_{l=-m}^{-1} \omega_q^{\,l\, \{ p \,(l+2m) + r \} } + \sum_{l=0}^{q-m-1} \omega_q^{\,l\, \{ p \,(l+2m) + r \} } \right].
\end{split}
\end{equation*}
 Now, it is easy to see that a shift of $q$ in $l$ leaves the sum unchanged, as every term only depends on the remainder of $l$ modulo $q$. Hence,
\begin{equation*}
    \begin{split}
        \sum_{l=-m}^{-1} \omega_q^{\,l\, \{ p \,(l+2m) + r \} } + \sum_{l=0}^{q-m-1} \omega_q^{\,l\, \{ p \,(l+2m) + r \} } &= \sum_{l=q-m}^{q-1} \omega_q^{\,l\, \{ p \,(l+2m) + r \} } + \sum_{l=0}^{q-m-1} \omega_q^{\,l\, \{ p \,(l+2m) + r \} }\\
        &= \sum_{l=0}^{q-1} \omega_q^{\,l\, \{ p \,(l+2m) + r \} }.
    \end{split}
\end{equation*}
Hence, we get
\begin{equation*}
    \abs{\mathbb{G}}^2 = \sum_{m=0}^{q-1} \sum_{l=0}^{q-1} \omega_q^{\,l\, \{ p \,(l+2m) + r \} } = \sum_{l=0}^{q-1} \left[ \omega_q^{\,pl^2 + rl } \sum_{m=0}^{q-1} \omega_q^{\,2pl\,m }\right] .
\end{equation*}
Now, the second sum is a geometric sum, whose general value is
\begin{equation*}
    \sum_{i=0}^{N} a\,b^i = a \,\dfrac{1-b^{N+1}}{1-b}.
\end{equation*}
In this case, $a=1$ and $b=\omega_q^{\, 2pl}$. Hence,
\begin{equation*}
    \abs{\mathbb{G}(p,r,q)}^2 = \sum_{l=0}^{q-1} \omega_q^{\,p l^2 + rl} \,\dfrac{1-\omega_q^{\, 2plq}}{1-\omega_q^{\, 2pl}} .
\end{equation*}
Now, the numerator of the fraction is always null, as $\omega_q^q=1$, thus the fraction will equal zero whenever the denominator is not null too. By inspection, it is easy to see that $1-\omega_q^{\, 2pl}=0$ is equivalent to $ 2pl \equiv 0\pmod{q}$, and as $(p,q)=1$, $2l \equiv 0\pmod{q}$. This only happens if $l=\{0; q/2\}$. In those cases, the geometric sum trivially equals $q$. Note that the case $l=q/2$ is only possible if $q$ is even.\smallskip

First, we will deal with the case where $q$ is odd, which is easier. Here we only consider that the fraction is not null for $l=0$, consequently,
\begin{equation*}
    \abs{\mathbb{G}(p,r,q)}^2 = q\,\omega_q^{\,p\,0^2 + r\,0} = q.
\end{equation*}
Hence,
\begin{equation}\label{caso1}
     q \equiv 1 \pmod{2} \implies \abs{\mathbb{G}}^2 = q.
\end{equation}
Otherwise, if $q$ is even, we must take into account also the term for $l=q/2$, so,
\begin{equation*}
    \abs{\mathbb{G}(p,r,q)}^2 = q\,\omega_q^{\,p\,0^2 + r\,0} + q\,\omega_q^{\,p\, (\frac{q}{2})^2 + r\frac{q}{2} } = q \left( 1 + e^{\frac{i\pi}{2} \, (p\, q + 2r)}\right)  = q \left( 1 + i^{\,p\, q + 2r }\right) .
\end{equation*}
Here, two cases arise. If $pq + 2r \equiv 0 \pmod{4} \implies i^{\,p\, q + 2r\,}=1 \implies \abs{\mathbb{G}}^2 = 2q$. But
\begin{equation*}
    pq + 2r \equiv 0 \pmod{4} \implies pq \equiv -2r \pmod{4} \implies \dfrac{pq}{2} \equiv -r \pmod{2}.
\end{equation*}
As $(p,q)=1$ and $2\mid q$, then $p$ is odd. So,
\begin{equation*}
    \dfrac{q}{2} \equiv r \pmod{2} \implies q \equiv 2r \pmod{4},
\end{equation*}
and, 
\begin{equation}\label{caso2}
    q \equiv 2r \pmod{4} \implies \abs{\mathbb{G}}^2 = 2q.
\end{equation}
Similarly, $pq + 2r \equiv 2 \pmod{4} \implies i^{\,p\, q + 2r\,}=-1 \implies \abs{\mathbb{G}}^2 = 0$. Thus,
\begin{equation}\label{caso3}
     q \not\equiv 2r \pmod{4} \implies \abs{\mathbb{G}}^2 = 0.
\end{equation}
Hence, by combining equations \eqref{caso1}, \eqref{caso2} and \eqref{caso3}, we get that,
\begin{empheq}[box=\fbox]{equation}
    \abs{\mathbb{G}(p,r,q)} = \begin{cases}
        \sqrt{q} 	& \text{if } q \text{ is odd},\\
        \sqrt{2q} 	& \text{if } 2 \mid q, \quad q \equiv r\pmod{4},\\
        0 			& \text{if } 2 \mid q, \quad q \not\equiv r\pmod{4}.
    \end{cases}
\end{empheq}
\hfill $\square$

%%%%%%%%%%%%%%%%%%%%%%%%% APENDICE SUMAS GAUSS MODIFICADAS

\newpage
\section{Proof of \texorpdfstring{$\abs{\mathbb{G}\left(p/2, pq/2 - m , q \right)} = \sqrt{q}$}{G = q}}\label{appendix:ideal}

Having defined
\begin{equation*}
    \mathbb{G}\left(p/2, pq/2 - m , q \right) := \sum_{r=0}^{q-1} e^{2\pi i \frac{(qp/2 + m) r - (p/2) r^2 }{q}},
\end{equation*}
we will prove that
\begin{empheq}[box=\fbox]{equation}
    \mathbb{G}\left(p/2, pq/2 - m , q \right) = \sqrt{q},\quad \forall p,q \in \mathbb{N}, \quad s.t. \gcd(p,q) = 1.
\end{empheq}
This proof is based on Appendices A and B of \cite{Eceizabarrena2021} and has been included in this work for completeness' sake. We will follow a case by case treatment, but first note that 
\begin{equation*}
    \mathbb{G}\left(p/2, pq/2 - m , q \right) = \sum_{r=0}^{q-1} e^{2\pi i (\frac{p}{2} + \frac{m}{q}) r - \pi i\frac{p}{q} r^2} = \sum_{r=0}^{q-1} e^{2\pi i (\frac{p \,(\text{mod}\,2)}{2} + \frac{m}{q}) r - \pi i\frac{p}{q} r^2}.
\end{equation*}
\subsection{The case \texorpdfstring{$p$}{p} even}
If $p \equiv 0\pmod{2}$, we may write 
\begin{equation*}
    \mathbb{G}\left(p/2, pq/2 - m , q \right) = \sum_{r=0}^{q-1} e^{2\pi i (\frac{m}{q}) r - \pi i\frac{p}{q} r^2}.
\end{equation*}
Thus $\abs{\mathbb{G}\left(p/2, pq/2 - m , q \right)} = \abs{\mathbb{G}(-p/2,m,q)}$, where $\mathbb{G}(p,r,q)$ is a generalized quadratic Gauss sum, whose absolute value is computed in Appendix \ref{appendix:GQGS}. The case $q$ even is impossible, as $(p,q)=1$. And if $q$ is odd, $\abs{\mathbb{G}(-p/2,m,q)} = \sqrt{q}$, as expected.

\subsection{The case \texorpdfstring{$p$}{p} odd}
By direct computation, 
\begin{equation*}
\begin{split}
    \mathbb{G}\left(p/2, pq/2 - m , q \right) &= \sum_{r=0}^{q-1} e^{2\pi i \frac{(q/2 + m)r - p/2 r^2}{q}} = \sum_{r=0}^{q-1} e^{2\pi i \frac{(q + 2m)r - p r^2}{2q}}\\
    &= \dfrac{1}{2}\, \sum_{r=0}^{2q-1} e^{2\pi i \frac{(q + 2m)r - p r^2}{2q}} = \dfrac{1}{2}\, \mathbb{G}(-p, 2m+q, 2q).
\end{split}
\end{equation*}

\subsubsection{The case \texorpdfstring{$q$}{q} odd}
Now, if $q$ is odd, we can use the multiplicative property $\mathbb{G}(a,b,cd) = \mathbb{G}(ac,b,d) \mathbb{G}(ad,b,c)$. So we may write
\begin{equation*}
    \abs{\mathbb{G}\left(p/2, pq/2 - m , q \right)} = \dfrac{1}{2}\abs{\mathbb{G}(-2p, 2m+q, q)} \abs{\mathbb{G}(-pq, 2m+q, 2)}
\end{equation*}
As $q$ is odd, $\abs{\mathbb{G}(-2p, 2m+q, q)}=\sqrt{q}$. Also, $pq$ and $2m+q$ must be odd, so\\
$\abs{\mathbb{G}(-pq, 2m+q, 2)} = \abs{\mathbb{G}(1, 1, 2)} = 2$. So $\abs{\mathbb{G}\left(p/2, pq/2 - m , q \right)} = \sqrt{q}$.

\subsubsection{The case \texorpdfstring{$q$}{q} even}
Finally, if $p$ odd and $q$ even, then, by completing the square,
\begin{equation*}
    \begin{split}
        \sum_{r=0}^{2q-1} e^{2\pi i\frac{(q+2m)r - p r^2}{2q}} &= \sum_{r=0}^{2q-1} e^{-2\pi i \frac{p \left(r -\frac{q+2m}{p} \right)^2}{2q}} e^{2\pi i \frac{p^{-1} (q+2m)^2}{2q}}\\
        &= e^{2\pi i \frac{p^{-1} (q+2m)^2}{2q}} \sum_{s=0}^{2q-1} e^{2\pi i \frac{-p s^2}{2q}} \\
        &= e^{2\pi i \frac{p^{-1} (q+2m)^2}{2q}} \, \mathbb{G}(-p, 0, 2q) .
    \end{split}
\end{equation*}
So $\abs{\mathbb{G}(-p,2m+q,2q)} = \abs{\mathbb{G}(-p,0,2q)} = 2\sqrt{q}$, as $\abs{\mathbb{G}(p,r,q)}=\sqrt{2q}$ if $q \equiv r\pmod{4}$. This finally gives us $\abs{\mathbb{G}\left(p/2, pq/2 - m , q \right)} = \sqrt{q}$. \hfill$\square$

\end{document}